\documentclass{siamart250211}
\usepackage{amsmath, amsfonts, amssymb, mathrsfs, bm}
\usepackage{mathtools}
\usepackage{microtype}
\usepackage{url}
\usepackage{graphicx}
\usepackage{booktabs}
\usepackage{placeins}
\usepackage{algorithm}
\usepackage{algorithmic}

\newsiamremark{remark}{Remark}
\AddToHook{env/remark/begin}{\crefalias{theorem}{remark}}

\def\ccm{Center for Computational Mathematics, Flatiron Institute, Simons Foundation, New York, New York 10010}
\def\courant{Courant Institute of Mathematical Sciences, New York University, New York, New York 10012}

\headers{An Adaptive Fast Algorithm for Periodic Lattice Sums}{X.~Gao, L.~Greengard, and S.~Jiang}

\newcommand{\ba}{\mathbf{a}}      % Lattice basis vector
\newcommand{\bk}{\mathbf{k}}      % Integer wave vector
\newcommand{\bn}{\mathbf{n}}      % Integer shift vector
\newcommand{\bx}{\mathbf{x}}      % Physical Target position
\newcommand{\by}{\mathbf{y}}      % Physical Source position
\newcommand{\bu}{\mathbf{u}}      % Lattice Target position
\newcommand{\bv}{\mathbf{v}}      % Lattice Source position
\newcommand{\bz}{\mathbf{z}}      % Difference vector

\newcommand{\bxi}{\bm{\xi}}       % Fourier variable (physical)
\newcommand{\bbeta}{\bm{\beta}}   % Bloch/Phase Vector
\newcommand{\qpot}{\Phi}          % Periodic potential (rendered \Phi; macro name kept)

\newcommand{\invL}{L^{-1}}        % Inverse Lattice Matrix
\newcommand{\detL}{|\det L|}      % Jacobian / Cell Volume
\newcommand{\rc}{r_c}             % Cutoff / cube side length
\newcommand{\bbu}{\bbeta_u}       % Bloch vector in fractional coords
\newcommand{\bmm}{\mathbf{m}}     % Integer Fourier mode vector
\newcommand{\C}{\mathcal{C}}      % Unit cell (matches pfmm)
\newcommand{\E}{\mathcal{E}}      % Extended domain
\newcommand{\pfmm}{\textsf{pfmm2d}} % pfmm2d software name
\newcommand{\npsi}{\tilde{\psi}_0^c}% Normalized prolate \psi_0^c/\psi_0^c(0)

\title{An Adaptive Fast Algorithm for Periodic Coulomb Lattice Sums in Arbitrary Unit Cells}

\author{
  Xuanzhao Gao\thanks{\ccm\,(\email{xgao@flatironinstitute.org}).}
  \and
  Leslie Greengard\thanks{\ccm;\,\courant\,(\email{lgreengard@flatironinstitute.org}).}
  \and
  Shidong Jiang\thanks{\ccm\,(\email{sjiang@flatironinstitute.org}).}
}

\begin{document}

\maketitle

\begin{abstract}
  We present a fast algorithm for evaluating conditionally
  convergent Coulomb lattice sums, governed by the Laplace equation
  with periodic boundary conditions on 
  arbitrary unit cells (oblique in 2D, triclinic in 3D) 
  and arbitrary particle distributions.  
  The algorithm extends the dual-space
  multilevel kernel-splitting (DMK) framework to this context.
  The root of the
  adaptive tree is now a rectangular grid of cubes consisting of an
  inner block covering the unit cell and a surrounding halo of image
  cubes, rather than a single cube, and the smooth top-level periodic
  kernel---the only term that requires the consideration of conditional
  convergence issues ---is
  evaluated by the ``five-step procedure" used in fast Ewald
  summation: spreading, fast Fourier transform (FFT), diagonal
  scaling, inverse FFT, and interpolation.  The resulting
  complexity is $O(N)$ for fixed cell shape.  
  Benchmarked against the periodic fast multipole method on
  highly nonuniform source distributions, our 2D algorithm is roughly
  an order of magnitude faster across particle counts and target
  precisions; in three dimensions, it is often as fast as the
  free-space DMK on the same sources, even for triclinic cells with
  edge-length ratios up to roughly $17$.
\end{abstract}



% ---------------------------------------------------------
\section{Introduction}
% ---------------------------------------------------------

Applications in electrostatics, molecular dynamics, and computational
materials science often involve point sources or volumetric
charge distributions contained in a unit
cell on which periodic boundary conditions are imposed. We consider
such a problem for the Laplace kernel in $d=2$ and $d=3$ dimensions,
restricting our attention to the discrete (point source) case.
Let the unit cell $\C \subset \mathbb{R}^d$ be the parallelepiped
(parallelogram in 2D) spanned by $d$ linearly independent lattice
vectors $\ba_1, \dots, \ba_d$, which we view as the columns of the lattice matrix
$L = [\ba_1, \dots, \ba_d] \in \mathbb{R}^{d \times d}$, and suppose that $N$
point sources are given, with strengths $q_j$ and positions
$\by_j \in \C$. The periodic potential at a target point
$\bx \in \mathbb{R}^d$ is given (formally) by the infinite lattice sum
\begin{equation}
    \qpot(\bx) = \sum_{\bn \in \mathbb{Z}^d} \sum_{j=1}^N G(\bx, \by_j + L\bn)\, q_j,
    \label{eq:potential}
\end{equation}
where $L\bn$ is the lattice translation associated with the integer
shift $\bn \in \mathbb{Z}^d$, and $G$ is the free-space Green's
function for the negative Laplacian:
\begin{equation}
    G(\bx, \by) =
    \begin{cases}
        -\dfrac{1}{2\pi}\,\log \|\bx - \by\|, & d=2, \\[4pt]
        \dfrac{1}{4\pi\,\|\bx - \by\|}, & d=3.
    \end{cases}
    \label{eq:greens}
\end{equation}
When the target $\bx$ coincides with one of the source positions
$\by_j$, the corresponding singular self-interaction term (with
$\bn = \mathbf{0}$) is excluded from the sum~\cref{eq:potential}. The
lattice sum is only \emph{conditionally convergent}. Throughout the
paper we assume both charge neutrality, $\sum_{j=1}^N q_j = 0$
and full periodicity of the potential. This corresponds to the 
usual macroscopic (tin-foil) convention that the 
$\bk = \mathbf{0}$ Fourier mode is zero and the solution is representable
by a Fourier series. In electrostatic applications, sums of the form
\eqref{eq:potential} are often called {\em Coulomb lattice sums}.

The challenge in evaluating~\cref{eq:potential} arises from the slow, conditional
convergence of the series.  Fast Ewald
methods~\cite{Darden1993,Essmann1995,hockney1988,Lindbo2011c,shamshirgar_fast_2021,
liang2026esp}
split the kernel into a smoothed far-field part that is evaluated by
the fast Fourier transform (FFT) on a uniform grid, and a compactly
supported near-field
residual that is summed directly.  They attain $O(N \log N)$
complexity when the sources are more or less uniformly distributed;
for highly clustered distributions, however, the global Ewald
parameter cannot be adjusted locally, the direct near-field sum
becomes the bottleneck, and the total cost can approach
$O(N^2)$~\cite{jiang2025cpam}.

Periodic boundary conditions are sometimes regarded as awkward to
impose within the fast multipole method (FMM).
However, for a cubic unit cell, the original
FMM~\cite{berman1994jmp,Greengard1987,greengard1987thesis} already addressed
that case in some detail.
For the kernel-independent
FMM~\cite{Ying2004}, singly, doubly, and triply periodic cubic geometries
can be handled 
via a precomputed multipole-to-local operator~\cite{Yan2018a}. For unit cells of arbitrary
shape, however, the situation is more delicate, and was
only recently addressed in~\cite{pei_fast_2023}. There, a free-space FMM is
augmented with low-rank compression of the field due to all
image sources beyond the first neighboring layer, derived from the
classical plane-wave (Sommerfeld integral) representation of the
free-space Green's function. The representation decays along a
single Cartesian direction, so that one such expansion is required per
coordinate direction.

In this work we extend the dual-space multilevel kernel-splitting
(DMK) framework~\cite{jiang2025cpam}---a fully adaptive $O(N)$
algorithm for free-space convolutions---to arbitrary unit cells~$\C$ in two
and three dimensions, while preserving its adaptivity to highly
nonuniform source distributions.  One main change is that the root
of the adaptive tree is now a
rectangular grid of cubes of side $\rc$ rather than a single cube; 
this multi-cube root contains an inner
block of $m_1 \times \cdots \times m_d$ cubes covering the unit cell
together with a surrounding halo of cubes carrying periodic copies
of the original sources, on which the near-field DMK is run.  The DMK
telescoping decomposition writes the kernel as a smooth far-field
component $W_0$ plus a sequence of difference and residual kernels
that are compactly supported at each tree level.  Since only $W_0$
has global support, the periodized $W_0^{(p)}$ is the only term
requiring a lattice sum: we represent it by a single Fourier
series on the reciprocal lattice, and evaluate its interaction by the five-step
procedure in~\cite{liang2026esp}---spreading, FFT, diagonal scaling,
inverse FFT, and interpolation. Equivalently, the procedure consists
of a type-1 nonuniform fast Fourier transform
(NUFFT)~\cite{nufft2,nufft3,Greengard2004}, 
multiplication by the Fourier transform of $W_0^{(p)}$ for each Fourier mode, and
a type-2 NUFFT, with both NUFFTs performed without upsampling. The compact support of the
difference and residual kernels then allows the remaining, near-field
contribution to be assembled as a \emph{free-space} DMK calculation
on the original sources together with this bounded image halo; the
ratio of image to original sources is bounded by a small geometric
constant depending on the cell shape and~$\rc$ (at most $16$ for the
parameter choices used in our experiments).
Recently, a related DMK-based approach was developed in \cite{krantz2026},
which permits periodicity in selected coordinate 
directions, but assumes an orthogonal unit cell
(of arbitrary aspect ratio).

Compared with the FMM-based scheme of~\cite{pei_fast_2023}---which,
like ours, is fully adaptive and handles highly nonuniform source
distributions---the present method has two structural advantages. First, the periodic
far-field requires only a \emph{single} Fourier series, in place of
$2d$ separate directional plane-wave expansions ($4$ in two
dimensions, $6$ in three), making the implementation noticeably
simpler. Second, the parameter $\rc$ which determines the 
cube size at the root level is a free parameter that
can be tuned to balance near-field and far-field work and to adapt
the method to the source nonuniformity; in particular, choosing
$\rc$ small enough that no root-level cube contains more than, say,
$n_s$ sources, would collapse the adaptive tree to a single
uniform grid and reduces the algorithm to fast Ewald summation. The present
scheme can be viewed as unifying the FMM and fast Ewald
summation (as the original DMK algorithm has done for the free-space case).

At fixed cell shape, the algorithm has $O(N)$ complexity.
Extensive numerical experiments on highly nonuniform 2D and 3D source
distributions show that our algorithm delivers $10$- to $26$-fold
speedups over the FMM-based scheme on 2D oblique cells with aspect
ratios up to $100$; in three
dimensions, it is often as fast as the underlying free-space DMK on
the same sources, even for triclinic cells with edge-length ratios
up to roughly $17$. To support reproducible evaluation, we release an open-source code
\texttt{PeriodicDMK}~\cite{PeriodicDMK}, written in C++ with a Julia
interface. Quasi-periodic
(Bloch) boundary conditions are supported with three localized
modifications, collected in \cref{sec:qp}.

The remainder of this paper is organized as follows.
\Cref{sec:prelim} sets up notation and recalls the prolate
spheroidal wave functions used as the windowing kernel together
with the dual-space multilevel kernel splitting
in~\cite{jiang2025cpam}.  \Cref{sec:pbc} extends the DMK
construction to periodic boundary conditions: the multi-cube root,
the image-source halo, the periodic far-field via the five-step
procedure of fast Ewald summation, and the modifications needed for
quasi-periodic (Bloch) boundary conditions.  \Cref{sec:algorithm}
states the full procedure as Algorithm~\ref{alg:pdmk} and
analyzes its computational and memory complexity, and
\cref{sec:choice_params} discusses the empirical tuning of the
algorithm's two free parameters.
\Cref{sec:numerics} reports numerical experiments on highly
nonuniform 2D and 3D source distributions and
~\cref{sec:conclusions} contains some concluding remarks.

% ---------------------------------------------------------
\section{Preliminaries}
\label{sec:prelim}
% ---------------------------------------------------------

\subsection{Notation and coordinate transformation}
\label{sec:coord}

Given lattice
vectors $\ba_1, \dots, \ba_d$, we define (as above) the lattice matrix
as $L = [\ba_1, \dots, \ba_d] \in \mathbb{R}^{d \times d}$.
It is easy to verify that $L$ maps the unit cell $\C$ onto the unit cube $[0,1)^d$.
The cell volume is denoted by  $V = \detL$.
For any target $\bx$ and source $\by$ in $\mathbb{R}^d$, we define the
\emph{fractional coordinates} by
\begin{equation}
    \bu = \invL \bx, \qquad \bv = \invL \by,
    \label{eq:fractional_coords}
\end{equation}
Under this mapping, Euclidean distance is
expressed via the metric tensor $M = L^T L$,
\begin{equation}
    \|\bx - \by\|^2 = (\bu - \bv)^T M (\bu - \bv).
    \label{eq:metric}
\end{equation}
This will permit us to process source/target interactions
on a standard Cartesian grid in fractional coordinates, as shown in
\cref{sec:farfield}.

Throughout the paper, the Fourier transform and its inverse are
defined by
\begin{equation}
    \widehat{f}(\bk) = \int_{\mathbb{R}^d} f(\bx)\,e^{-i\bk\cdot\bx}\,d\bx,
    \qquad
    f(\bx) = \frac{1}{(2\pi)^d}\int_{\mathbb{R}^d} \widehat{f}(\bk)\,e^{i\bk\cdot\bx}\,d\bk.
    \label{eq:ft_convention}
\end{equation}

\subsection{Prolate spheroidal wave functions}
\label{sec:pswf}

Fix a bandwidth parameter $c > 0$. The zeroth prolate spheroidal
wave function (PSWF) $\psi_0^c$ is the eigenfunction corresponding
to the largest eigenvalue $\lambda_0$ of the integral operator
$\varphi \mapsto \int_{-1}^{1} \varphi(t)\,e^{i c x t}\,dt$ on
$L^2([-1,1])$~\cite{landau1961bstj,slepian1961bstj}; that is,
\begin{equation}
    \int_{-1}^{1} \psi_0^c(t)\,e^{i c x t}\,dt
    = \lambda_0\,\psi_0^c(x),
    \qquad x \in [-1,1].
    \label{eq:pswf_eig}
\end{equation}
$\psi_0^c$ is even and positive on $(-1,1)$, and 
we adopt the convention that it is extended to $\mathbb{R}$ by zero outside
$[-1,1]$. Its Fourier transform is given by 
\begin{equation}
    \widehat{\psi}_0^{\,c}(k)
    = \lambda_0\,\psi_0^c(k/c),
    \qquad |k| \le c,
    \label{eq:pswf_ft}
\end{equation}
so $\widehat{\psi}_0^{\,c}$ is also essentially supported on $[-c,c]$.
Among bandlimited functions of bandwidth $c$, $\psi_0^c$ is optimal
in the sense that it maximizes the fraction of $L^2$ energy concentrated
in $[-1,1]$. We calibrate the bandwidth $c$ to a target precision
$\varepsilon$ by $\psi_0^c(1) \approx \varepsilon$, which yields the
scaling
\begin{equation}
    c \;\sim\; \log(1/\varepsilon).
    \label{eq:c_log_eps}
\end{equation}
For convenience we introduce the normalized prolate
\begin{equation}
    \npsi(x) := \frac{\psi_0^c(x)}{\psi_0^c(0)},
    \label{eq:npsi}
\end{equation}
so that $\npsi(0) = 1$.

\subsection{Dual-space multilevel kernel splitting (DMK)}
\label{sec:dmk}

We briefly recall the telescoping kernel decomposition from the DMK
framework~\cite{jiang2025cpam}. A defining feature of DMK is that
the splitting can be carried out in either physical or Fourier
space; for the Laplace kernel, the Fourier-space form is most
convenient, since the Fourier transform of the Green's function~\cref{eq:greens} is
\begin{equation}
  \widehat{G}(\bxi) = \frac{1}{|\bxi|^2} \, ,
  \label{eq:Ghat}
\end{equation}
independent of dimension.

Let the box size of the cubes at the root level be denoted by
$h_0$. Superimposing a level-restricted adaptive quad-tree (2D) or oct-tree (3D)
data structure on the cubes, the box size at level $\ell$ is given by 
$h_\ell = h_0/2^\ell$. 

With $\npsi$ the normalized prolate defined in
\cref{eq:npsi}, the DMK splitting of $\widehat{G}$ is
\begin{equation}
  \widehat{G}(\bxi)
  = \widehat{W}_0(\bxi)
  + \sum_{\ell=0}^{\mathscr{L}-1} \widehat{D}_\ell(\bxi)
  + \widehat{R}_{\mathscr{L}}(\bxi),
  \label{eq:tele}
\end{equation}
with components
\begin{equation}
  \begin{aligned}
    \widehat{W}_0(\bxi) &= \frac{\npsi(|\bxi|\,h_0/c)}{|\bxi|^2}, \\
    \widehat{D}_\ell(\bxi) &= \frac{\npsi(|\bxi|\,h_{\ell+1}/c) - \npsi(|\bxi|\,h_\ell/c)}{|\bxi|^2}, \\
    \widehat{R}_{\mathscr{L}}(\bxi) &= \frac{1 - \npsi(|\bxi|\,h_{\mathscr{L}}/c)}{|\bxi|^2}.
  \end{aligned}
  \label{eq:W0hat}
\end{equation}
Since $\psi_0^c$ extends to an even entire function on
$\mathbb{C}$ that is essentially supported on $[-1,1]$, each
$\widehat{D}_\ell$ is smooth at $\bxi = \mathbf{0}$ 
and essentially supported on $|\bxi| \le c/h_\ell$.
(The apparent $1/|\bxi|^2$ singularity is canceled since the numerator is zero
at the origin and even.)  
Thus, from the Paley-Wiener theorem, 
the difference kernel $D_\ell$ is exponentially small in
physical space once $|\bx-\by| > h_\ell$. The same holds for
the residual kernel $R_{\mathscr{L}}$ on $|\bx-\by| \le h_{\mathscr{L}}$.

For direct interactions evaluated in physical space (involving
$R_{\mathscr{L}}$), the residual kernel $R_{\mathscr{L}}(r)$ is required pointwise. In three
dimensions, the inverse Fourier transform of \cref{eq:W0hat} admits
a closed form: 
\begin{equation}
  W_0(r) = \frac{\phi_0^c(r/h_0)}{4 \pi r}, \quad
  D_\ell(r) = \frac{\phi_0^c(r/h_{\ell+1}) - \phi_0^c(r/h_\ell)}{4 \pi r}, \quad
  R_{\mathscr{L}}(r) = \frac{1 - \phi_0^c(r/h_{\mathscr{L}})}{4 \pi r},
  \label{eq:pswf_split}
\end{equation}
where
\begin{equation}
  \phi_0^c(t) = \frac{1}{c_0}\int_0^t \psi_0^c(s)\,ds,
  \qquad c_0 = \int_0^1 \psi_0^c(s)\,ds.
  \label{eq:Phic}
\end{equation}
We note that $W_0$ 
is smooth everywhere (including at the origin) and satisfies the property $W_0(r) = G(r)$ for $r\ge h_0$.  
In two dimensions the logarithmic singularity of $G$ at the origin
prevents an analogous closed form, and we instead obtain $R_{\mathscr{L}}(r)$ by
numerical inverse Fourier transform of $\widehat{R}_{\mathscr{L}}$, using the
radial Fourier transform approach of~\cite{jiang2025cpam}.

In DMK, at each
level $\ell$ of the adaptive tree, the \emph{colleagues} of a box are the same-level boxes
sharing a boundary point. (There are up to $3^d$ such boxes in $d$ dimensions). By
construction, $D_\ell$ interactions are confined to colleagues at
level $\ell$, $R_{\mathscr{L}}$ interactions to colleagues at the leaf level, and
$W_0$ accounts for interactions beyond the colleague set at the root
level. 

% ---------------------------------------------------------
\section{DMK for periodic boundary conditions}
\label{sec:pbc}
% ---------------------------------------------------------

We now describe how to extend the DMK framework to compute the
periodic potential $\qpot(\bx)$ in~\cref{eq:potential} for an
arbitrary unit cell defined by $L$ in $d$ dimensions.  The key
insight is that the compactly supported components of the splitting
(the difference kernels $D_\ell$ and the residual $R_{\mathscr{L}}$) only
require nearby periodic images, while the smooth far-field component
$W_0$ is efficiently represented by a Fourier series on the
reciprocal lattice.

\begin{figure}[htbp]
  \centering
  \includegraphics[width=0.8\linewidth]{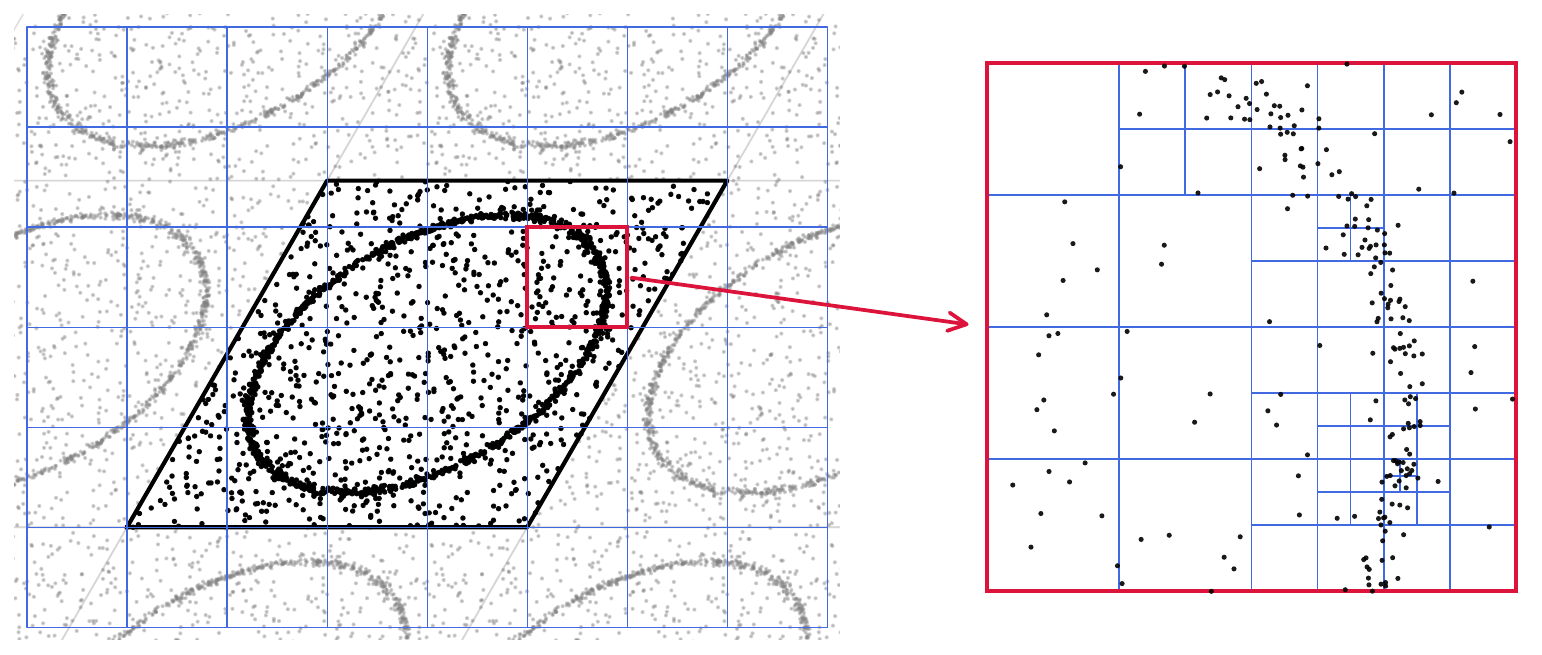}
  \caption{Schematic overview of the algorithm in $d=2$ for an
    oblique unit cell with interior angle $\pi/3$ and a highly
    nonuniform source distribution.
    \textbf{Left:} the unit cell (shaded parallelogram) is embedded
    in the smallest enclosing rectangular grid of $6 \times 4$
    cubes of side~$\rc$, then surrounded by a one-cube halo of
    periodic image sources, yielding the extended
    $(6+2)\times(4+2)$ root grid on which the near-field DMK is
    run.
    \textbf{Right:} the level-restricted (2:1 balanced) adaptive
    quadtree built within a single root cube, illustrating how
    refinement follows the local source density and how the
    balancing is enforced across root-cube boundaries.}
  \label{fig:overview}
\end{figure}

\subsection{Periodized kernel splitting}

For any kernel $F$ defined on $\mathbb{R}^d$, denote by $F^{(p)}$ its
formal periodization,
\begin{equation}
  F^{(p)}(\bx,\by)
  = \sum_{\bn \in \mathbb{Z}^d} F(\bx - \by - L\bn).
  \label{eq:periodize}
\end{equation}
Applying periodization term by term to the multilevel kernel
splitting of the Green's function in~\cref{sec:dmk} gives
\begin{equation}
  G^{(p)}(\bx,\by) = W_0^{(p)}(\bx,\by)
  + \sum_{\ell=0}^{\mathscr{L}-1} D_\ell^{(p)}(\bx,\by)
  + R_{\mathscr{L}}^{(p)}(\bx,\by).
  \label{eq:pbc_split}
\end{equation}
The compactly supported components $D_\ell$ and $R_{\mathscr{L}}$ have
absolutely convergent periodizations, with only finitely many shifts
$\bn$ contributing for each $(\bx,\by)$; only the smooth far-field
kernel $W_0$ retains the conditional-convergence issue, which is
resolved in \cref{sec:farfield} via Poisson summation.

The potential therefore can be expressed as
\begin{equation}
  \qpot(\bx_i)
  = \underbrace{\sum_{j=1}^N W_0^{(p)}(\bx_i,\by_j)\, q_j}
    _{\displaystyle\qpot_{\mathrm{far}}(\bx_i)}
  \;+\;
  \underbrace{\sum_{j=1}^N \left[
    \sum_{\ell=0}^{\mathscr{L}-1} D_\ell^{(p)}(\bx_i,\by_j)
    + R_{\mathscr{L}}^{(p)}(\bx_i,\by_j)\right] q_j}
    _{\displaystyle\qpot_{\mathrm{near}}(\bx_i)}.
  \label{eq:near_far}
\end{equation}
We describe the evaluation of each contribution in turn.

\subsection{Embedding the unit cell in a rectangular grid}
\label{sec:embed}

The construction below is presented for $d = 3$; the $d = 2$ case is
entirely analogous (with one fewer coordinate and an upper-triangular
$2\times 2$ lattice matrix).

Without loss of generality, we orient the
lattice so that $\ba_1$ lies along the positive $x_1$-axis with
$|\ba_1| \ge |\ba_2| \ge |\ba_3|$, and $\ba_2$ lies in the upper
half of the $x_1 x_2$-plane ($x_2 \ge 0$).  Any user-supplied basis
can be brought to this form by first permuting the lattice vectors
in order of decreasing magnitude, then applying a single orthogonal
transformation $Q \in O(d)$---a rigid rotation when the input basis
is right-handed, or a rotation composed with a reflection when it is
left-handed.  Since $Q$ preserves all pairwise distances, the
potential $\qpot$ is invariant.  The resulting
lattice matrix is upper triangular with positive diagonal:
\begin{equation}
  L = \begin{pmatrix}
    L_{11} & L_{12} & L_{13} \\
    0      & L_{22} & L_{23} \\
    0      & 0      & L_{33}
  \end{pmatrix},
  \qquad L_{11} \ge L_{22} > 0, \quad L_{33} > 0,
  \label{eq:Ltri}
\end{equation}
where $L_{11} = |\ba_1|$.

\paragraph{Bounding box of the unit cell}
The unit cell $\C = \{L\bu : \bu \in [0,1)^3\}$ has one corner at
the origin.  Since $L$ is upper triangular, a point
$\bx = L\bu$ has coordinates
\begin{equation}
  x_1 = L_{11}\,u_1 + L_{12}\,u_2 + L_{13}\,u_3, \quad
  x_2 = L_{22}\,u_2 + L_{23}\,u_3, \quad
  x_3 = L_{33}\,u_3.
  \label{eq:coords}
\end{equation}
We compute the range of each coordinate over $u_1, u_2, u_3 \in
[0,1)$, working from the bottom row upward.

\emph{Third coordinate.}
$x_3 = L_{33}\,u_3$ with $L_{33} > 0$ and $u_3 \in [0,1)$, so
\begin{equation}
  x_3^{-} = 0, \qquad x_3^{+} = L_{33}.
  \label{eq:x3range}
\end{equation}

\emph{Second coordinate.}
$x_2 = L_{22}\,u_2 + L_{23}\,u_3$ is a sum of two independent terms.
Since $L_{22} > 0$, the term $L_{22}\,u_2$ ranges over $[0, L_{22})$.
The term $L_{23}\,u_3$ ranges over $[0, L_{23})$ if $L_{23} \ge 0$,
or over $(L_{23}, 0]$ if $L_{23} < 0$.  Adding the two ranges gives
\begin{equation}
  x_2^{-} = \min(0,\,L_{23}), \qquad
  x_2^{+} = L_{22} + \max(0,\,L_{23}).
  \label{eq:x2range}
\end{equation}
The width in the $x_2$-direction is
\begin{equation}
  w_2 = x_2^{+} - x_2^{-}
  = L_{22} + \max(0,\,L_{23}) - \min(0,\,L_{23})
  = L_{22} + |L_{23}|,
  \label{eq:w2}
\end{equation}
using the identity $\max(0,a) - \min(0,a) = |a|$.

\emph{First coordinate.}
$x_1 = L_{11}\,u_1 + L_{12}\,u_2 + L_{13}\,u_3$ is a sum of three
independent terms.  Since $L_{11} > 0$, the first ranges over
$[0, L_{11})$.  The second and third contribute $[\min(0,L_{12}),
\max(0,L_{12}))$ and $[\min(0,L_{13}), \max(0,L_{13}))$ respectively.
Adding gives
\begin{equation}
  x_1^{-} = \min(0,\,L_{12}) + \min(0,\,L_{13}), \qquad
  x_1^{+} = L_{11} + \max(0,\,L_{12}) + \max(0,\,L_{13}),
  \label{eq:x1range}
\end{equation}
with width
\begin{equation}
  w_1 = x_1^{+} - x_1^{-}
  = L_{11} + |L_{12}| + |L_{13}|.
  \label{eq:w1}
\end{equation}
We collect the per-coordinate extrema into the corner vectors
$\bx^{-} = (x_1^{-}, x_2^{-}, x_3^{-})$ and
$\bx^{+} = (x_1^{+}, x_2^{+}, x_3^{+})$, so that the bounding box
of the unit cell is $[\bx^{-}, \bx^{+}]$.

\paragraph{Covering with cubes}
Let $\rc$ be a given cutoff radius (chosen so that the difference
kernel $D_0$ is supported within a ball of radius $\rc$); in the
notation of \cref{sec:dmk}, $\rc$ plays the role of the
root-level box size~$h_0$.  We cover the bounding box
$[\bx^{-}, \bx^{+}]$ by a grid of cubes of side length~$\rc$.  The grid can be placed at
an arbitrary offset: its faces in the $x_i$-direction lie at
$s_i + k\,\rc$ for $k \in \mathbb{Z}$, where $s_i$ is a free shift
parameter.  Choosing $s_i = x_i^{-}$ aligns one grid face with the
left edge of the bounding box, and the number of cubes needed in
direction~$i$ is simply $\lceil w_i / \rc \rceil$:
\begin{equation}
  m_3 = \left\lceil \frac{L_{33}}{\rc} \right\rceil, \qquad
  m_2 = \left\lceil \frac{L_{22} + |L_{23}|}{\rc} \right\rceil,
    \qquad
  m_1 = \left\lceil \frac{L_{11} + |L_{12}| + |L_{13}|}{\rc}
    \right\rceil.
  \label{eq:mi}
\end{equation}
The resulting grid of $m_1 \times m_2 \times m_3$ cubes of side $\rc$
is the smallest such grid that encloses the unit cell.

The cutoff $\rc$ is a free parameter, subject to the upper bound
$\rc \le \min_i w_i$ that ensures $m_i \ge 1$ for all $i$. The
simplest (and safest) default is
\begin{equation}
  \rc = \min_i L_{ii} \equiv L_{\min},
  \label{eq:rc_choice}
\end{equation}
which yields $m_i = \lceil L_{ii}/L_{\min}\rceil$ for orthogonal
cells, directly reflecting the aspect ratio. However, this choice
is rarely optimal in practice: a substantially smaller $\rc$
typically improves both runtime and peak memory. The choice of
$\rc$ together with the leaf capacity $n_s$ is discussed
in detail in~\cref{sec:choice_params}.

\subsection{The multi-cube root and extended domain}
\label{sec:multicube}

The DMK framework extends transparently to the case where the root
of the adaptive tree consists of a rectangular grid of cubes rather
than a single cube.  Starting from the cover of the unit cell~$\C$
by $m_1 \times \cdots \times m_d$ cubes constructed in
\cref{sec:embed}, we surround it by one additional layer of cubes in
each direction, yielding the \emph{extended domain} $\E$ of
\begin{equation}
  (m_1 + 2) \times (m_2 + 2) \times \cdots \times (m_d + 2)
  \text{ cubes of side } \rc,
  \label{eq:extended}
\end{equation}
on which the periodic near-field will be evaluated.  The halo is
populated with \emph{image sources}: for each lattice shift
$\bn \in \mathbb{Z}^d$ such that $\by_j + L\bn \in \E$, a copy of
source~$j$ is placed at $\by_j + L\bn$ with charge
\begin{equation}
  q_j^{(\bn)} = q_j.
  \label{eq:image_charge}
\end{equation}
The enumeration of contributing shifts is detailed in
\cref{sec:images}.

With the extended grid as the root, the DMK tree refines, sorts
particles, and runs its upward and downward passes exactly as in the
free-space case; only the level-$0$ data structures differ.  Each
root cube has up to $3^d$ colleagues (the box itself together with
its same-level neighbors) among the
$(m_1+2) \times \cdots \times (m_d+2)$ root cubes (subject to the
rectangular boundary), and the level-restricted ($2{:}1$) balance
must be enforced across root-cube boundaries rather than within a
single root cube; see \cref{alg:pdmk} for the algorithmic details.

\subsection{Image-source enumeration}
\label{sec:images}

The image sources fill all of $\E$ outside the unit cell~$\C$: a
one-cube halo around the $m_1\times\cdots\times m_d$ inner block,
plus---for triclinic cells, where $\C$ does not fill its bounding
box---all interior cubes of the bounding box that lie outside~$\C$.
We enumerate the contributing lattice shifts $\bn \in \mathbb{Z}^d$,
defined as those producing at least one translated source
$\by_j + L\bn$ ($\by_j \in \C$) inside~$\E$.

Since $L$ is upper triangular, the $i$-th component of $L\bn$ is
$\sum_{k \ge i} L_{ik}\,n_k$, so any image landing in $\E$
satisfies $|\sum_{k \ge i} L_{ik}\,n_k| \lesssim (m_i+1)\,\rc$.
We enumerate shifts in nested loops with $n_d$ outermost; after
fixing the outer indices, the innermost is bounded by
\begin{equation}
  n_i^{\max}(n_{i+1}, \ldots, n_d)
  \;=\; \left\lceil \frac{(m_i + 1)\,\rc
      + \bigl\lvert \sum_{k=i+1}^{d} L_{ik}\,n_k \bigr\rvert}
      {L_{ii}} \right\rceil,
  \qquad i = d, d-1, \ldots, 1,
  \label{eq:shift_bounds}
\end{equation}
with the empty-sum convention giving
$n_d^{\max} = \lceil (m_d+1)\rc/L_{dd} \rceil$.  Any shift violating
$|n_i| \le n_i^{\max}$ is pruned without per-source work; for each
shift $\bn$ that survives the bound, every source~$\by_j$ is then
subjected to the \emph{in-box test}: accept $\by_j + L\bn$ as an
image if and only if it lies in~$\E$.

The total source count in~$\E$---originals together with all
accepted images---is bounded by $N$ times a distribution-independent
constant determined by the cell geometry and $\rc$
(\cref{eq:Nbound}); a looser bound on $n_i^{\max}$ simply prunes
more shifts at the in-box test rather than producing more accepted
images.  The cubic specialization is treated separately
in~\cref{rem:cube}.

\begin{remark}[Cube case]
\label{rem:cube}
The cubic-cell specialization has already been treated
in~\cite{afklinteberg2026jcp}, where the cutoff $\rc$ was chosen
as a quarter of the cube side, so that the periodized windowed
kernel $W_0^{(p)}$ absorbs all far-field interactions at the first
two levels of a single-cube-root DMK tree.  In this
setting there is no need for the multi-cube root construction of
\cref{sec:multicube}, and no image sources need to be added
to the tree or sorted: a standard single-cube DMK tree is built
on the $N$ original particles alone, and image interactions at
each boundary box are handled by introducing virtual neighbors
whose positions are shifted by $L\bn$ with
$\bn \in \{-1,0,1\}^d$, instead of allocating and sorting physical
image particles.  This eliminates
the image-source array, the extra tree-build pass, and the
reordering cost altogether.  For a general triclinic cell, by
contrast, the off-diagonal entries of~$L$ displace image positions
orthogonally to~$\bn$, so that images no longer coincide with
grid points at the coarsest tree levels, and the full image
construction of this subsection is required.
\end{remark}

\subsection{Near-field: free-space DMK on the extended domain}
\label{sec:nearfield}

Since every difference kernel $D_\ell$ is supported within a ball of
radius at most $\rc$ (with progressively smaller support at deeper
levels), and the residual kernel $R_{\mathscr{L}}$ has even smaller support, all
near-field interactions $\qpot_{\mathrm{near}}(\bx_i)$
in~\cref{eq:near_far} for targets $\bx_i \in \C$ involve only
sources within distance~$\rc$.  By construction, every such
source---whether original or image---lies in~$\E$.

Consequently, the evaluation of $\qpot_{\mathrm{near}}$ reduces to a
\emph{standard free-space DMK calculation} on~$\E$, treating the
original and image sources uniformly.  Using the
multi-cube root from \cref{sec:multicube}, we build a DMK tree
with root consisting of
$(m_1+2) \times \cdots \times (m_d+2)$ cubes, and run the standard
upward pass, downward pass (for difference kernels $D_\ell$), and
direct evaluation (for the residual kernel $R_{\mathscr{L}}$).  Only the
\emph{root-level far-field step} (the $W_0$ convolution) is
\emph{omitted}, as it is replaced by the periodic Fourier series
described next.

To verify correctness, note that for any target $\bx_i\in\C$
and source~$\by_j$, the periodized near-field involves the sum
\[
  \sum_{\bn \in \mathbb{Z}^d}
  \bigl[D_\ell(\bx_i - \by_j - L\bn) + R_{\mathscr{L}}(\bx_i - \by_j - L\bn)
  \bigr].
\]
Terms with $\|\bx_i - \by_j - L\bn\| > \rc$ vanish by compact
support.  All remaining terms correspond to image sources inside~$\E$.
Thus the free-space DMK on~$\E$ reproduces the exact periodized
near-field.

\cref{fig:overview} summarizes the geometric scaffolding implied by
the splitting~\cref{eq:pbc_split}.  The left panel shows the extended
root grid: the unit cell~$\C$ embedded in $m_1 \times m_2$ cubes of
side~$\rc$ surrounded by a one-cube image halo, on which the
compactly supported $D_\ell^{(p)}$ and $R_{\mathscr{L}}^{(p)}$ contributions
reduce to a free-space DMK call.  The right panel shows the
level-restricted adaptive tree inside one root cube, with $2{:}1$
balance enforced even across root-cube boundaries.  The smooth
far-field $W_0^{(p)}$ has no geometric footprint in this picture:
it is evaluated entirely in reciprocal space by a pair of NUFFTs,
described next.

\subsection{Far-field: Fourier series and the five-step procedure}
\label{sec:farfield}

The periodized windowed kernel $W_0^{(p)}$ admits an absolutely
convergent Fourier series on the reciprocal lattice.  Define the
reciprocal lattice vectors $\bk_\bmm = 2\pi L^{-T}\bmm$ for
$\bmm \in \mathbb{Z}^d$.  Applied to
$W_0^{(p)}(\bz) = \sum_{\bn} W_0(|\bz - L\bn|)$ with
$\bz = \bx - \by$, the Poisson summation formula yields
\begin{equation}
  W_0^{(p)}(\bx,\by)
  = \frac{1}{V}\sum_{\bmm \in \mathbb{Z}^d \setminus \{\mathbf{0}\}}
    \widehat{W}_0(|\bk_\bmm|)\;
    e^{i\bk_\bmm\cdot(\bx - \by)},
  \label{eq:W0p}
\end{equation}
where $V = \detL$ is the cell volume and the $\bmm = \mathbf{0}$ term
is omitted (see ``Charge neutrality'' below).  Since
$\widehat{W}_0(|\bxi|)$ from~\cref{eq:W0hat} is essentially supported
on $|\bxi| \le c/h_0 = c/\rc$, the series~\cref{eq:W0p} is
effectively finite: only modes with $|\bk_\bmm| \lesssim c/\rc$
contribute to precision~$\varepsilon$.  The corresponding rectangular
NUFFT grid that contains all active modes is fixed in
\cref{sec:modes}.

The far-field potential
$\qpot_{\mathrm{far}}(\bx_i) = \sum_{j=1}^N W_0^{(p)}(\bx_i,\by_j)\,q_j$
is evaluated by the standard \emph{five-step procedure} of fast Ewald
summation~\cite{liang2026esp}---spreading, FFT, diagonal scaling,
inverse FFT, and interpolation.  Using the fractional coordinates of
\cref{sec:coord} (in which $\bk_\bmm\cdot\bx = 2\pi\,\bmm\cdot\bu$),
spreading and FFT bundle into a type-1
NUFFT~\cite{nufft2,nufft3,Greengard2004}, and inverse FFT and
interpolation bundle into a type-2 NUFFT, giving the equivalent
three-step realization:
\begin{enumerate}
  \item \textbf{Type-1 NUFFT.}  From the charges $q_j$ at the
    non-uniform fractional source positions
    $\bv_j = \invL\by_j \in [0,1)^d$, compute
    \[
      \widehat{g}(\bmm) = \sum_{j=1}^N q_j\,
      e^{-2\pi i\,\bmm \cdot \bv_j}
    \]
    for all active modes $\bmm$ (spreading to a uniform grid followed
    by an FFT).

  \item \textbf{Diagonal scaling.}  Multiply each Fourier coefficient
    by the kernel weight,
    \begin{equation}
      \widehat{g}(\bmm) \;\leftarrow\;
      \frac{1}{V}\;\widehat{W}_0(|\bk_\bmm|)\;
      \widehat{g}(\bmm),
      \label{eq:diag}
    \end{equation}
    omitting the $\bmm = \mathbf{0}$ mode.

  \item \textbf{Type-2 NUFFT.}  Evaluate
    \[
      \qpot_{\mathrm{far}}(\bx_i)
      = \sum_{\bmm} \widehat{g}(\bmm)\,
      e^{2\pi i\,\bmm \cdot \bu_i}
    \]
    at the non-uniform fractional target positions
    $\bu_i = \invL\bx_i$ (an inverse FFT on the uniform grid followed
    by interpolation).
\end{enumerate}
For real charges and real targets the output is real.

\paragraph{Charge neutrality}
The $\bmm = \mathbf{0}$ mode involves $\widehat{W}_0(0) \to \infty$.
Its exclusion is justified when $\sum_{j=1}^N q_j = 0$ (charge
neutrality), which is a necessary condition for the conditionally
convergent periodic Coulomb lattice sum to be well-defined.

\begin{remark}[No upsampling]
\label{rem:no_upsampling}
The standard NUFFT uses an oversampled grid (typically by a factor
of~2 in each dimension) to suppress aliasing in the
spreading/interpolation stages.  In the present setting this
oversampling is unnecessary: the diagonal scaling~\cref{eq:diag}
multiplies $\widehat{g}(\bmm)$ by $\widehat{W}_0(|\bk_\bmm|)/V$, which
decays spectrally fast beyond $|\bk_\bmm| \gtrsim c/\rc$, so the
Fourier coefficients after scaling are already negligible outside the
active mode set.  Both NUFFTs can therefore operate on a grid of
exactly $M_1 \times \cdots \times M_d$ points (\cref{sec:modes}),
saving a factor of~$2^d$ in the FFT cost.  This observation is
discussed in detail in~\cite{liang2026esp} in the context of Ewald
summation with prolate spheroidal wave functions.
\end{remark}

\subsection{Choosing the Fourier grid}
\label{sec:modes}

The NUFFT operates on a rectangular grid of integer mode vectors
$\bmm = (m_1, \ldots, m_d)$ with $m_i$ ranging from $-M_i/2$ to
$(M_i-1)/2$.  In $\bmm$-space the active-mode constraint
$|\bk_\bmm| \lesssim c/\rc$ from~\cref{sec:farfield} defines an
ellipsoid, since $\bk_\bmm = 2\pi L^{-T}\bmm$, and the grid
dimensions $M_1, \ldots, M_d$ must be large enough that the
rectangular box contains this ellipsoid.

\paragraph{Derivation}
Inverting $\bk_\bmm = 2\pi L^{-T}\bmm$ gives
$\bmm = (2\pi)^{-1} L^T \bk_\bmm$, whose $i$-th component is
\begin{equation}
  m_i = \frac{1}{2\pi}\,\ba_i \cdot \bk_\bmm,
  \label{eq:mi_dot}
\end{equation}
since $\ba_i$ is the $i$-th column of $L$ and the $i$-th row of $L^T$
equals $\ba_i^T$.  Applying the Cauchy--Schwarz inequality to the
bandwidth constraint $|\bk_\bmm| \le c/\rc$ yields
\begin{equation}
  |m_i|
  \le \frac{|\ba_i|\,|\bk_\bmm|}{2\pi}
  \le \frac{c\,|\ba_i|}{2\pi\,\rc},
  \label{eq:mi_bound}
\end{equation}
with equality when $\bk_\bmm$ is parallel to~$\ba_i$.  The NUFFT grid
dimensions are therefore
\begin{equation}
  M_i = 2\left\lceil \frac{c\,|\ba_i|}{2\pi\,\rc} \right\rceil,
  \qquad i = 1, \ldots, d,
  \label{eq:Mi}
\end{equation}
the smallest even integer grid that contains all active modes in
direction~$i$.  By \cref{rem:no_upsampling} this is also the grid
size used by the NUFFT itself, with no further oversampling.

\paragraph{Lattice vector lengths}
Since $L$ is upper triangular, the $i$-th lattice vector
$\ba_i$ (column of $L$) has length
$|\ba_i| = \bigl(\sum_{k=1}^{i} L_{ki}^2\bigr)^{1/2}$.
In $d = 3$ this gives explicitly
\begin{equation}
  |\ba_1| = L_{11}, \qquad
  |\ba_2| = \sqrt{L_{12}^2 + L_{22}^2}, \qquad
  |\ba_3| = \sqrt{L_{13}^2 + L_{23}^2 + L_{33}^2},
  \label{eq:anorms}
\end{equation}
and in $d = 2$,
$|\ba_1| = L_{11}$, $|\ba_2| = \sqrt{L_{12}^2 + L_{22}^2}$.
For an orthogonal cell (all off-diagonal $L_{ij} = 0$), this reduces
to $|\ba_i| = L_{ii}$ and the Fourier grid is isotropic in the
natural sense: $M_i = 2\lceil c\,L_{ii}/(2\pi\,\rc)\rceil$.
With the choice $\rc = L_{\min} = \min_i L_{ii}$ from
\cref{eq:rc_choice}, the grid dimension in the $i$-th direction is
$M_i \approx c\,L_{ii}/(\pi\,L_{\min})$, directly reflecting the
aspect ratio of the cell.

\paragraph{Overhead for skewed cells}
The NUFFT grid is rectangular, while the active modes occupy an
ellipsoid.  The ratio of their volumes is
\begin{equation}
  \frac{\prod_i M_i}{N_F}
  \;\approx\; \frac{2^d}{S_d}\;
  \frac{\prod_i|\ba_i|}{V},
  \label{eq:mode_overhead}
\end{equation}
where $S_d$ is the volume of the unit ball in $\mathbb{R}^d$
($S_2 = \pi$, $S_3 = 4\pi/3$).
By Hadamard's inequality,
$\prod_i|\ba_i| \ge V = \det L$, with equality if and only
if the lattice vectors are mutually orthogonal.  For orthogonal cells,
the overhead factor is $4/\pi \approx 1.27$ in $d = 2$ and
$6/\pi \approx 1.91$ in $d = 3$; for moderately skewed
cells it remains modest.  For highly skewed cells, the overhead can
grow, but this is an inherent cost of using a rectangular Fourier grid
on a non-orthogonal lattice.

\subsection{Self-interaction correction}
\label{sec:selfcorr}

The near-field and far-field contributions of
\cref{sec:nearfield,sec:farfield} assemble to the periodic
potential
\begin{equation}
  \qpot(\bx_i) = \qpot_{\mathrm{near}}(\bx_i)
  + \qpot_{\mathrm{far}}(\bx_i),
  \label{eq:assembly}
\end{equation}
after one final correction.  The plane-wave interactions in the DMK
hierarchy include self-interaction terms: at each level, the smooth
kernel does not vanish at zero separation, and the standard DMK
subtracts these self-terms at each leaf.  In the PBC setting, the
NUFFT evaluation of $W_0^{(p)}$ implicitly adds a further spurious
$W_0(\mathbf{0})\,q_j$ to the potential at every target $\bx_i$
coinciding with a source $\by_j$.

The value $W_0(\mathbf{0})$ has a closed form in $d=3$.  Using
$W_0(r) = \phi_0^c(r/\rc)\,G(r)$ from~\cref{eq:pswf_split} together
with the small-$t$ expansion $\phi_0^c(t) \approx \psi_0^c(0)\,t/c_0$,
\begin{equation}
  W_0(\mathbf{0})
  = \frac{\psi_0^c(0)}{c_0\,\rc}\;\lim_{r\to 0}\,r\,G(r)
  = \frac{\psi_0^c(0)}{4\pi\,c_0\,\rc}.
  \label{eq:W0zero}
\end{equation}
In $d = 2$, $W_0$ is obtained by numerical inverse Fourier transform
of $\widehat{W}_0$, and $W_0(\mathbf{0})$ is read off by the same
procedure (see~\cite{jiang2025cpam} for details).  The leaf-level
self-subtraction already in DMK is then extended to subtract
$W_0(\mathbf{0})\,q_j$ at each coincident source--target pair as
well, after which the assembly~\cref{eq:assembly} returns the
correct periodic potential.

% ---------------------------------------------------------
\section{Algorithm and complexity}
\label{sec:algorithm}
% ---------------------------------------------------------

\Cref{alg:pdmk} summarizes the complete procedure.  Derivations
of each step appear in~\cref{sec:pbc}, and the modifications
needed for quasi-periodic boundary conditions are collected
in~\cref{sec:qp}.

\begin{algorithm}[ht]
\caption{Periodic DMK}\label{alg:pdmk}
\begin{algorithmic}[1]
\REQUIRE Lattice matrix $L$, sources $\{(\by_j, q_j)\}_{j=1}^N$,
targets $\{\bx_i\}_{i=1}^M$, target precision $\varepsilon$.
\STATE \emph{Preprocessing.}  Orient $L$ as in~\cref{sec:embed} so
that it is upper-triangular (\cref{eq:Ltri}).  Choose
$\rc = L_{\min}/\kappa$ and the leaf capacity $n_s$
from~\cref{tab:tuning}.  Compute the bounding-box widths $w_i$ and
root-grid dimensions $m_i$ from~\cref{eq:mi}.  Populate the
extended domain~$\E$ with image charges
via~\cref{eq:image_charge}.
\STATE \emph{Tree construction.}  Build the adaptive tree on~$\E$
with $(m_1+2)\times\cdots\times(m_d+2)$ root-level cubes of
side~$\rc$, recursively subdividing any box holding more than
$n_s$ sources into $2^d$ children of side $\rc/2$.  Patch the
level-$0$ colleague lists (\cref{sec:multicube}), then enforce
$2{:}1$ balance across root-cube
boundaries~\cite{sundar2008sisc}.
\STATE \emph{Near-field DMK.}  Run the free-space DMK on~$\E$
(\cref{sec:nearfield}) with the root-level windowed kernel $W_0$
omitted; by~\cref{eq:tele} the result is exactly $G - W_0$, the
near-field part of~\cref{eq:near_far}.
\STATE \emph{Far-field NUFFT.}  Apply the five-step procedure
of~\cref{sec:farfield} to the original sources at fractional
coordinates $\bv_j = \invL\by_j$ and targets at
$\bu_i = \invL\bx_i$, with the $\bmm = \mathbf{0}$ mode of the
diagonal scaling~\cref{eq:diag} omitted.
\STATE \emph{Assembly and self-correction.}  Form
$\qpot(\bx_i) = \qpot_{\mathrm{near}}(\bx_i)
+ \qpot_{\mathrm{far}}(\bx_i)$ as in~\cref{eq:assembly}, and
subtract $W_0(\mathbf{0})\,q_j$ at every coincident source--target
pair (\cref{sec:selfcorr}).
\RETURN $\qpot(\bx_i)$ for $i = 1, \ldots, M$.
\end{algorithmic}
\end{algorithm}

\subsection{Quasi-periodic boundary conditions}
\label{sec:qp}

A \emph{quasi-periodic} (Bloch) potential satisfies
\begin{equation}
  \qpot(\bx + L\bn) = e^{i\bbeta\cdot L\bn}\,\qpot(\bx),
  \qquad \bn \in \mathbb{Z}^d,
  \label{eq:qp_def}
\end{equation}
for a Bloch vector $\bbeta \in \mathbb{R}^d$; the periodic case is
recovered by setting $\bbeta = \mathbf{0}$.  The periodic construction
of~\cref{sec:pbc} extends to this setting through three localized
modifications.

\paragraph{Periodized kernel}
The periodized Green's function and its components acquire the
Bloch phase: $G^{(p)}(\bx,\by) = \sum_{\bn} G(\bx,\by + L\bn)\,
e^{i\bbeta\cdot L\bn}$, and likewise for $W_0^{(p)}$,
$D_\ell^{(p)}$, $R_{\mathscr{L}}^{(p)}$.  In the near-field DMK on~$\E$, this is
implemented by replacing the image-source charges in
\cref{eq:image_charge} with $q_j^{(\bn)} = q_j\,e^{i\bbeta\cdot L\bn}$;
no other near-field changes are required.  The far-field Fourier
series~\cref{eq:W0p} becomes
\begin{equation}
  W_0^{(p)}(\bx,\by)
  = \frac{1}{V}\sum_{\bmm \in \mathbb{Z}^d}
    \widehat{W}_0(|\bk_\bmm + \bbeta|)\,
    e^{i(\bk_\bmm + \bbeta)\cdot(\bx - \by)},
  \label{eq:W0p_qp}
\end{equation}
with all modes retained, provided $\bk_\bmm + \bbeta \ne \mathbf{0}$
for every $\bmm \in \mathbb{Z}^d$---equivalently, $\bbeta \notin
2\pi L^{-T}\mathbb{Z}^d$; this excludes only the Bloch vectors
equivalent to the periodic case modulo the reciprocal lattice.
Under this condition every $\widehat{W}_0(|\bk_\bmm + \bbeta|)$ is
finite, and charge neutrality is no longer required.

\paragraph{Fourier grid}
The active-mode condition becomes
$|\bk_\bmm + \bbeta| \lesssim c/\rc$, raising the grid
dimensions~\cref{eq:Mi} to
$M_i = 2\lceil(c/\rc + |\bbeta|)\,|\ba_i|/(2\pi)\rceil$.

\paragraph{Five-step procedure}
Let $\bbu = L^T\bbeta$ denote the Bloch vector in fractional
coordinates.  Three of the five steps of~\cref{sec:farfield} are
modified:
\begin{enumerate}
  \item \emph{Source pre-modulation.}  Replace the charges in the
    type-1 NUFFT with $\widetilde{q}_j = q_j\,e^{-i\bbu\cdot\bv_j}$.
  \item \emph{Diagonal scaling at shifted modes.}  Replace
    $\widehat{W}_0(|\bk_\bmm|)$ in~\cref{eq:diag} by
    $\widehat{W}_0(|\bk_\bmm + \bbeta|)$, and retain the
    $\bmm = \mathbf{0}$ mode.
  \item \emph{Target post-modulation.}  Multiply the type-2 output
    by $e^{i\bbu\cdot\bu_i}$.
\end{enumerate}
Combining these three steps with the image-source phasing above, the
full quasi-periodic far-field reads
\begin{equation}
  \qpot_{\mathrm{far}}(\bx_i)
  = \frac{e^{i\bbu\cdot\bu_i}}{V}
    \sum_{\bmm \in \mathbb{Z}^d}
    \widehat{W}_0(|\bk_\bmm + \bbeta|)\,
    \biggl(\sum_{j=1}^N q_j\,
      e^{-i\bbu\cdot\bv_j}
      e^{-2\pi i\,\bmm\cdot\bv_j}\biggr)
    e^{2\pi i\,\bmm\cdot\bu_i}.
  \label{eq:qp_far}
\end{equation}
No other algorithmic changes are needed.

\subsection{Complexity analysis}
\label{sec:complexity}

We take the number of targets equal to the number of sources,
$M = N$; the case $M \ne N$ replaces $N$ by $N + M$ in the
near-field DMK passes, while the image-halo bound below depends
only on the source count~$N$.  The computational cost has two
components.

\paragraph{Near-field (free-space DMK)}
The near-field DMK on~$\E$ operates on $N'$ total sources (original
plus images).  Each particle $\by_j \in \C$ contributes a copy at
$\by_j + L\bn$ for every shift $\bn \in \mathbb{Z}^d$ with $\by_j +
L\bn \in \E$, and the number of admissible shifts is bounded above
by a constant determined by the cell geometry and $\rc$, independent
of the source distribution.  For an orthogonal cell with edges
$L_i$ and $m_i = \lceil L_i/\rc \rceil$ cubes along axis $i$, an
elementary count along each axis gives
\begin{equation}
  \frac{N'}{N} \;\le\;
  \prod_{i=1}^{d}\!\left(
    \left\lfloor (m_i + 2)\,\tfrac{\rc}{L_i} \right\rfloor + 1
  \right);
  \label{eq:Nbound}
\end{equation}
analogous axis-by-axis bounds for triclinic cells follow from the
image-shift enumeration of~\cref{eq:shift_bounds}.  The bound
\cref{eq:Nbound} is distribution-independent even though $N'$
itself is not, so the near-field cost is $O(N)$ with $N'/N$ bounded
by a geometric constant that depends on the cell shape and~$\rc$.
For the tuning parameter $\rc = L_{\min}/\kappa$ with $\kappa \ge 1$
in our experiments (\cref{tab:tuning}), the bound evaluates to a
small number --- at most $16$ in the worst case (the $d = 3$ slab
regime with $\kappa = 1$).

\paragraph{Far-field (NUFFT)}
Let $N_F$ denote the number of active Fourier modes, i.e., those
$\bmm \in \mathbb{Z}^d$ with $|2\pi L^{-T}\bmm| \lesssim c/\rc$.
Counting lattice points in the active ball gives
\begin{equation}
  N_F \approx \frac{(c/\rc)^d\,V}{(2\pi)^d}.
  \label{eq:NF}
\end{equation}
Define the aspect ratio
\begin{equation}
  A \;=\; L_{\max}/L_{\min},
  \label{eq:aspect_ratio}
\end{equation}
where $L_{\max}, L_{\min}$ are the largest and smallest diagonal
entries of~\cref{eq:Ltri}.  For an orthogonal cell with edges
$L_1 \ge \cdots \ge L_d$, $V/L_{\min}^d = \prod_i(L_i/L_{\min})$,
which lies between $A$ (a rod with one long edge) and $A^{d-1}$
(a $(d{-}1)$-dimensional slab).  Each NUFFT costs
$O(N_F \log N_F + N \log^d(1/\varepsilon))$.

\paragraph{Total cost}
Combining the near-field and far-field costs, the total complexity is
\[
  O(N) \;+\; O(N_F \log N_F),
\]
where $N_F$ is given by~\cref{eq:NF}.  With $\rc \propto L_{\min}$,
reciprocal-lattice mode counting gives $N_F = \Theta(V/L_{\min}^d)$,
a quantity fixed by the cell geometry and independent of~$N$.  The
algorithm is therefore $O(N)$, with the $O(N_F \log N_F)$ NUFFT term
contributing only a constant, cell-shape-dependent overhead.

Quantitatively, for an orthogonal cell with aspect ratio
$A = L_{\max}/L_{\min}$, $V/L_{\min}^d$ ranges from $A$ in the rod
case (one long edge) to $A^{d-1}$ in the slab case ($d{-}1$ long
edges), as discussed after~\cref{eq:NF}.  The fixed NUFFT overhead
is therefore $O(1)$ for the isotropic regime ($A = O(1)$), $O(A \log A)$
for the rod regime, and $O(A^2 \log A)$ for the slab regime in $d = 3$.
The triclinic cell carries an additional Hadamard factor
$\prod_i|\ba_i|/V$ from~\cref{eq:mode_overhead}, which still depends
only on cell shape.

Whether this overhead is negligible depends on $N_F$ relative
to~$N$: for $N_F \lesssim N$, the DMK upward and downward passes
dominate (\cref{tab:2d,tab:rect3d,tab:tri3d}); when $N_F \gg N$ ---
which requires very large $A$ in $d = 2$ or $A^2$ in $d = 3$ --- the
NUFFT overhead can dominate instead.

\paragraph{Memory}
The main memory consumers are the NUFFT coefficient grid of
$\prod_i M_i$ complex entries, the outgoing plane-wave expansions
on the tree, the tensor-product proxy grid on each leaf, and the
reordered source/charge arrays of length~$N'$.  The plane-wave
buffer can be reused across tree levels, because each cross-level
plane-wave translation reads only the current level's outgoing data
and writes into the next level's incoming buffer.  The concatenated
original-and-image source array assembled during preprocessing can
be released as soon as the tree-sort copies it into the working
buffers; it is not referenced afterward.  With these
observations, the peak working memory is bounded by the largest
single level of plane-wave data plus the NUFFT grid plus $O(N')$.

% ---------------------------------------------------------
\section{Choice of parameters}
\label{sec:choice_params}
% ---------------------------------------------------------

The two adjustable parameters of the algorithm are the root-cube
side length $\rc$ and the leaf capacity $n_s$.  A direct cost
analysis for general nonuniform source distributions is difficult,
since the per-level work depends on the local density at each scale
and varies from one distribution to another.  We therefore tune
$\rc$ and $n_s$ empirically --- a standard practice for FMM-type
algorithms.

\paragraph{Choice of $\rc$}
The cutoff $\rc$ divides the cost between the NUFFT far-field and
the near-field DMK on the extended domain~$\E$, and its optimum
depends sensitively on the cell geometry.  As $\rc$ decreases,
the active-mode count $N_F \approx (c/\rc)^d V/(2\pi)^d$
from~\cref{eq:NF} scales as $\rc^{-d}$, raising both the FFT work
and the NUFFT coefficient grid storage; the number of root-level
boxes grows in proportion, increasing the per-box computation and
the plane-wave coefficient storage.  As $\rc$ increases, the image halo grows in volume according to
$V_{\mathrm{ext}}/V \le (1 + 3\rc/L_{\min})^d$ for orthogonal
cells, raising the source count entering the near-field DMK; at
$\rc = L_{\max}$ the halo volume exceeds the cell volume by at
least a factor of $3^d$, with the factor growing further for cells
of high aspect ratio.  The total cost
is therefore $U$-shaped in $\rc$, and in practice we find that
$\rc = L_{\min}/\kappa$ for some $\kappa > 1$ works well across a
wide range of particle distributions and cell shapes.

\paragraph{Choice of $n_s$}
\label{sec:choice_ns}
The leaf capacity $n_s$ is an internal DMK parameter --- it controls
only the near-field tree and is unaffected by the periodic far-field
or the cell geometry, so its optimum is essentially that of free-space
DMK.  As $n_s$ increases, the direct-evaluation work at the leaves grows
linearly, while the plane-wave bookkeeping at the internal tree
levels shrinks because the tree becomes shallower with fewer boxes.  The cost is again
$U$-shaped, with optimum depending on the precision, the dimension,
and machine constants.

\begin{table}[!th]
  \centering
  \caption{Tuning parameters for the 2D and 3D experiments.
    $\varepsilon$: target relative $L^2$ error.
    $n_s$ is the leaf capacity; $\kappa$ sets the root-cube divisor
    $\rc = L_{\min}/\kappa$.  In 2D a single $\kappa{=}7$ is used
    across all $\varepsilon$ and aspect ratios $A\in\{1,10,100\}$.
    In 3D $\kappa$ varies with the anisotropy regime: isotropic
    (a)~$1{:}1{:}1$, rod (b)~$1{:}1{:}10$, and slab
    (c)~$1{:}10{:}10$.  Values were chosen by an offline sweep
    over a single representative configuration per (regime,
    $\varepsilon$).}
  \label{tab:tuning}

  \begin{tabular}{lcccccc}
    \toprule
    $\varepsilon$ & \multicolumn{2}{c}{2D} & \multicolumn{4}{c}{3D} \\
    \cmidrule(lr){2-3} \cmidrule(lr){4-7}
                   & $\kappa$ & $n_s$ & $\kappa$ (iso) & $\kappa$ (rod) & $\kappa$ (slab) & $n_s$ \\
    \midrule
    $10^{-3}$  & 7 & 30 & 7 & 5 & 1 & 80 \\
    $10^{-6}$  & 7 & 40 & 7 & 5 & 1 & 200 \\
    $10^{-9}$  & 7 & 40 & 7 & 3 & 1 & 500 \\
    $10^{-12}$ & 7 & 60 & 7 & 3 & 1 & 600 \\
    \bottomrule
  \end{tabular}
\end{table}

\paragraph{Empirical values}
The values used in~\cref{sec:numerics} are summarized
in~\cref{tab:tuning}; since $\kappa$ depends sensitively on the cell
geometry, they are specific to the cell shapes tested.  In every 3D
regime $\kappa$ is no larger than the 2D value: in 2D the
active-mode count scales as $\rc^{-2}$, so a small $\rc$ (large
$\kappa$) remains affordable, while in 3D it scales as $\rc^{-3}$
and the NUFFT/proxy-grid memory rapidly becomes limiting, favoring
a larger $\rc$ (smaller $\kappa$).  In 3D $\kappa$ shrinks further
as the cell becomes more anisotropic, since the periodic neighbor
count along the long axis grows like $L_{\max}/\rc$ and must be
held bounded.  Conversely, $n_s$ is larger in 3D than in 2D because
the residual kernel $R_{\mathscr{L}}$ is cheaper to evaluate per pair
in 3D, so a larger $n_s$ is tolerable before direct work dominates.

\begin{remark}[Adaptive tuning]
A more refined design would introduce a runtime planning stage
analogous to FFTW's~\cite{Frigo2005} that probes the source
distribution and the machine to pick $\rc$ and $n_s$ adaptively
per problem; we defer this to future work.  Like other multilevel fast algorithms, the
present scheme is memory-bound on modern hardware, so the empirical
optima reflect cache behavior and bandwidth as much as FLOP balance.
\end{remark}

% ---------------------------------------------------------
\section{Numerical results}
\label{sec:numerics}
% ---------------------------------------------------------

Our implementation, hereafter referred to as
PDMK~\cite{PeriodicDMK}, is built on top of the open-source
free-space DMK code of~\cite{DMKcode}, extended with the multi-cube
root, image-source construction, and periodic far-field evaluation
described in~\cref{sec:pbc} and summarized in~\cref{alg:pdmk}.  The type-1 and
type-2 NUFFTs used in the far-field step are provided by the FINUFFT
library~\cite{barnett2019sisc,finufftlib}.  All timings are single-thread, double
precision, on an AMD EPYC 9474F CPU with 4.1 GHz clock speed and
1.5 TB of memory.  Our implementation is built with
\texttt{-march=x86-64-v4}; the free-space DMK
baseline~\cite{DMKcode} and the periodic FMM in 2D
(\pfmm{}~\cite{pei_fast_2023}) are built with the same
compiler flags from their reference implementations.

We do not benchmark the cube case: it admits the specialized treatment
described in \cref{rem:cube} and has already been studied in that
setting.  We also do not benchmark uniform or nearly uniform particle
distributions: in that regime, single-level fast Ewald summation is
generally much faster than any adaptive scheme, including the present
algorithm.  The experiments below therefore focus on the regime where
the proposed method is designed to excel---\emph{highly nonuniform
particle distributions on non-cubic unit cells}.

In two dimensions, the unit cell is an oblique parallelogram
with interior angle $\pi/3$, i.e.\ lattice vectors
\[
  \ba_1 = (1,\,0), \qquad
  \ba_2 = a\,(\cos\tfrac{\pi}{3},\,\sin\tfrac{\pi}{3})
        = (a/2,\, a\sqrt{3}/2),
\]
where the aspect ratio $a = |\ba_2|/|\ba_1|$ takes the three fixed
values $1$, $10$, and $100$.  For each aspect ratio, we
sweep the particle count $N$ from $10^5$ to $10^7$.

In three dimensions, we use two \emph{base} cells with the same set of edge lengths $(1,\,\sqrt{2},\,\sqrt{3})$, which cannot be scaled to have an integer ratio.
The \emph{rectangular} cell has all three lattice angles
equal to $\pi/2$; the \emph{triclinic} cell has all three angles
equal to $\pi/3$.  Each base cell is then scaled by an integer
triple $(s_1, s_2, s_3)$, giving three regimes (the labels $\mathrm{(a)},
\mathrm{(b)}, \mathrm{(c)}$ match the rows of \cref{tab:rect3d,tab:tri3d}):
\begin{enumerate}
  \item[(a)] \emph{isotropic} $(1,1,1)$: the original cell, with
    edge lengths in the ratio $\sqrt{3}\,{:}\,\sqrt{2}\,{:}\,1$;
  \item[(b)] \emph{rod} $(1,1,10)$: the third edge is stretched by
    $10\times$, so the longest edge is about $17\times$ the
    shortest;
  \item[(c)] \emph{slab} $(1,10,10)$: the second and third edges
    are stretched by $10\times$, so the two longer edges are about
    $14\times$ and $17\times$ the shortest.
\end{enumerate}
The two cell types together with three scaling modes give six 3D test
geometries.

\begin{figure}[!th]
  \centering
  \includegraphics[width=0.8\linewidth]{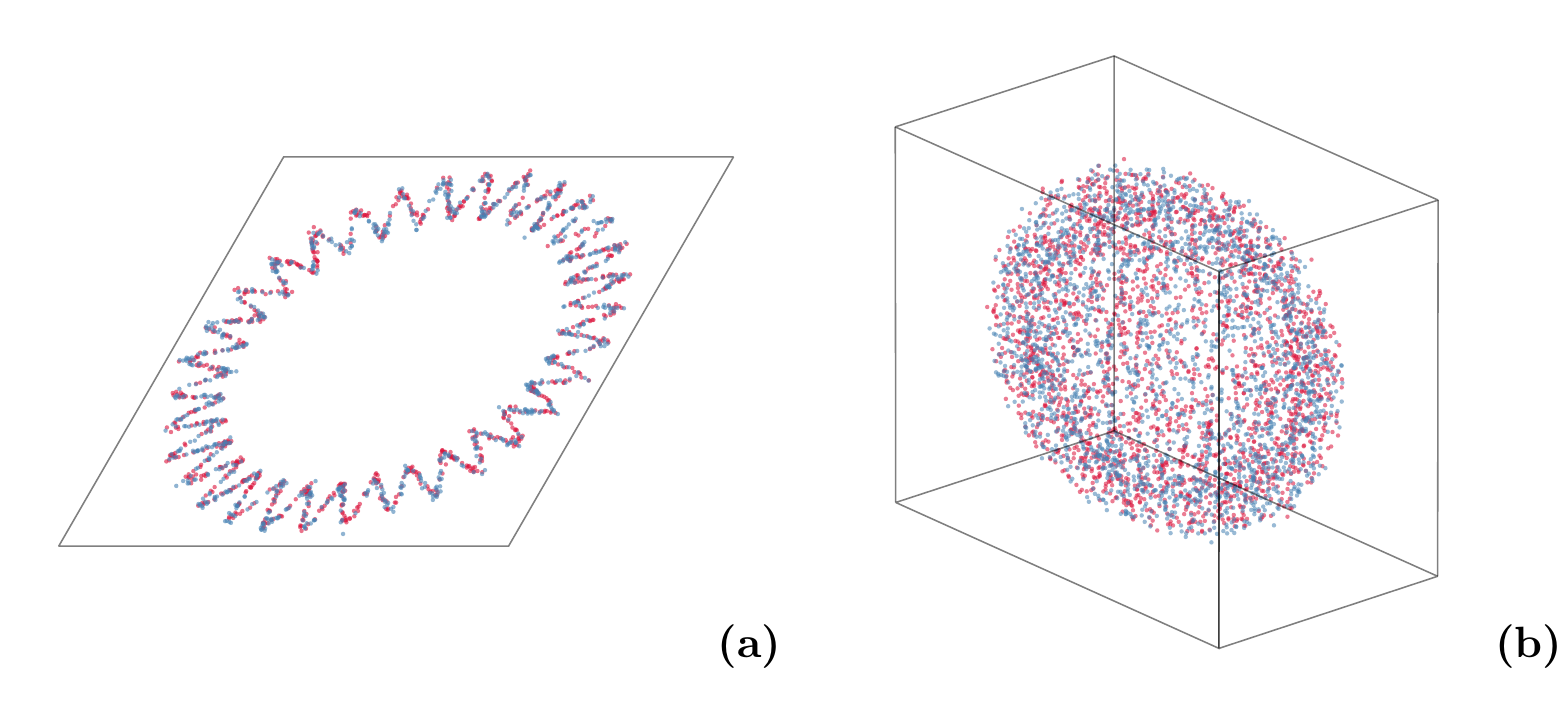}
  \caption{Representative 2D wavy-ring test system in the
    oblique ($\pi/3$) unit cell and 3D wavy-sphere system in the rectangular cell.
    For 2D, sources sit on a $32$-lobed
    perturbed ring concentric with the cell (drawn in black).
    The 3D wavy-sphere distribution is the natural azimuthal lift
    of this 2D system, notice that all particles are distributed near the surface instead of the interior.}
  \label{fig:systems}
\end{figure}

\paragraph{Source distributions}
Both test families place sources on a thin azimuthally-modulated
ring (2D) or spherical shell (3D) centered in the unit cell, with
fractional radius
$r_j = r_0 + A\sin(k\theta_j) + \eta_j$, where $\theta_j$ is
the azimuth, that is drawn uniformly on $[0, 2\pi)$ in 2D, or read off
from a uniformly-drawn direction on $S^2$ in 3D,
$\eta_j \sim \mathcal{N}(0, \sigma^2)$ is a small Gaussian radial
jitter, and the parameters $r_0 = 0.4$, $A = 0.05$, $k = 32$,
$\sigma = 0.01$ are fixed independently of $N$ and the cell aspect
ratio.  Charges are zero-mean, as required by charge neutrality
for the conditionally convergent periodic Coulomb lattice sum.  The $k=32$
modulation gives 32 narrow azimuthal lobes of angular width
$\approx 0.2$~rad and hence a strongly multi-scale source set whose
local density varies by close to an order of magnitude over very
short scales---precisely the regime where adaptive refinement pays
off and uniform-grid Ewald methods struggle.  \cref{fig:systems}
illustrates a representative 2D wavy-ring test system and the 3D
wavy-sphere distribution.

\paragraph{Reported quantities}
For each cell we report three sets of measurements at four target
precisions $\varepsilon \in \{10^{-3}, 10^{-6}, 10^{-9}, 10^{-12}\}$:
(i) wall-clock time versus particle count $N$ over the range
$10^5$ to $10^7$, with a dashed linear-fit reference;
(ii) average throughput in particles per second; and (iii) a stage-level timing
table splitting the total into $t_{\mathrm{build}}$,
$t_{\mathrm{nufft}}$, and $t_{\mathrm{eval}}$, with the periodic
FMM (2D) or free-space DMK (3D) listed alongside.  The components
are defined in each table caption.

\subsection{Two dimensions}
\label{sec:results2d}

\begin{figure}[!th]
  \centering
  \includegraphics[width=0.9\linewidth]{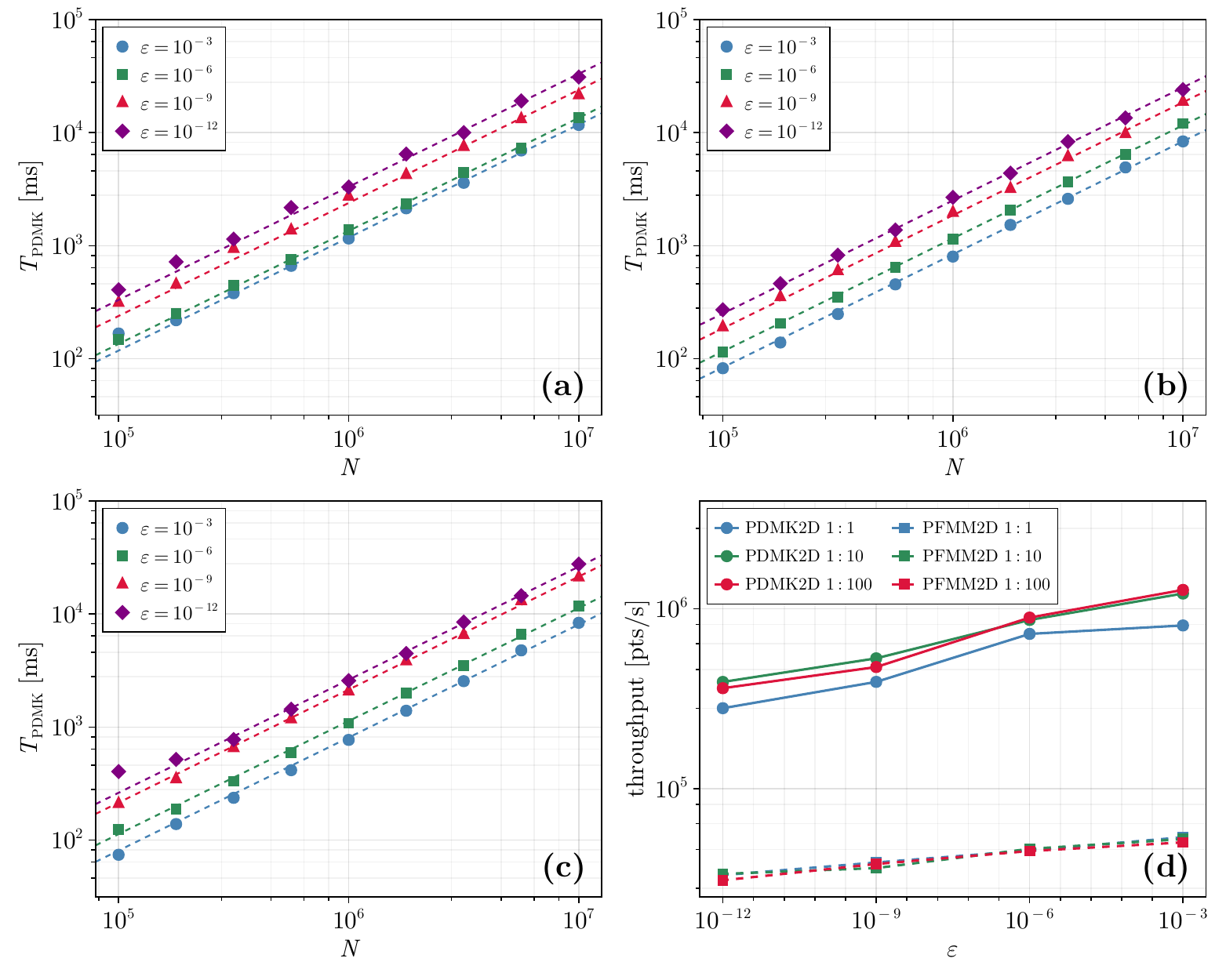}
  \caption{2D results on oblique ($\pi/3$) wavy-ring
    cells of increasing aspect ratio.  (a): runtime versus $N$ at
    four target precisions $\varepsilon = 10^{-3}, 10^{-6}, 10^{-9},
    10^{-12}$, with linear-fit reference lines, for aspect ratio $A = 1$.
    (b): same for $A = 10$.  (c): same for $A = 100$.
    (d): average throughput versus precision and aspect ratio of PDMK and \pfmm{}.}
  \label{fig:rt2d}
\end{figure}

\begin{table}[!th]
  \centering
  \caption{Per-step PDMK runtime and head-to-head comparison
    with \pfmm{}~\cite{pei_fast_2023} on the oblique
    ($\pi/3$) wavy-ring cells of \cref{fig:rt2d}.
    $N$: number of sources (equal to the number of targets) fixed to $10^7$.
    $A$: aspect ratio of the unit cell.
    $\varepsilon$: target precision.
    $t_{\mathrm{build}}$: tree construction with image-source
    placement.
    $t_{\mathrm{nufft}}$: type-1 NUFFT, diagonal scaling, and
    type-2 NUFFT for $W_0^{(p)}$.
    $t_{\mathrm{eval}}$: near-field DMK (plane-wave shifts,
    upward/downward passes, and leaf-level direct evaluation
    of $R_{\mathscr{L}}$).
    $t_{\mathrm{PDMK}} = t_{\mathrm{build}} + t_{\mathrm{nufft}}
    + t_{\mathrm{eval}}$: total PDMK time.
    $t_{\mathrm{PFMM}}$: total \pfmm{} time on the same
    system.
    ratio: PDMK throughput divided by \pfmm{} throughput,
    which (since both run on the same $N$) equals
    $t_{\mathrm{PFMM}}/t_{\mathrm{PDMK}}$; values above $1$ mean
    PDMK has higher throughput.
    $E_{\text{rel}}$: relative $L^2$ error.
    Times in seconds; single thread, double precision.}
  \label{tab:2d}
  \resizebox{\textwidth}{!}{\begin{tabular}{rlrrrrrrr}
  \toprule
  $A$ & $\varepsilon$
      & $t_{\mathrm{build}}$ & $t_{\mathrm{nufft}}$ & $t_{\mathrm{eval}}$ & $t_{\mathrm{PDMK}}$
      & $t_{\mathrm{PFMM}}$ & ratio
      & $E_{\mathrm{rel}}$ \\
  \midrule
  $1$ & $10^{-3}$ & 3.84 & 0.89 & 6.72 & 11.45 & 202.12 & 17.65 & $1.5\times 10^{-3}$ \\
  $1$ & $10^{-6}$ & 3.87 & 1.09 & 9.51 & 14.47 & 261.53 & 18.07 & $4.7\times 10^{-6}$ \\
  $1$ & $10^{-9}$ & 4.61 & 1.48 & 15.42 & 21.51 & 285.24 & 13.26 & $1.5\times 10^{-9}$ \\
  $1$ & $10^{-12}$ & 5.08 & 1.81 & 24.18 & 31.07 & 311.04 & 10.01 & $4.6\times 10^{-13}$ \\
  \midrule
  $10$ & $10^{-3}$ & 2.13 & 0.90 & 5.34 & 8.37 & 210.29 & 25.13 & $1.7\times 10^{-4}$ \\
  $10$ & $10^{-6}$ & 2.22 & 1.14 & 8.68 & 12.05 & 229.91 & 19.08 & $7.5\times 10^{-7}$ \\
  $10$ & $10^{-9}$ & 2.75 & 1.50 & 14.59 & 18.84 & 281.03 & 14.92 & $2.7\times 10^{-10}$ \\
  $10$ & $10^{-12}$ & 2.51 & 1.88 & 21.61 & 26.00 & 335.79 & 12.91 & $9.1\times 10^{-14}$ \\
  \midrule
  $100$ & $10^{-3}$ & 1.95 & 0.97 & 5.45 & 8.37 & 219.11 & 26.17 & $1.3\times 10^{-4}$ \\
  $100$ & $10^{-6}$ & 2.14 & 1.20 & 8.42 & 11.76 & 255.35 & 21.71 & $1.2\times 10^{-7}$ \\
  $100$ & $10^{-9}$ & 2.82 & 1.61 & 15.19 & 19.62 & 298.67 & 15.22 & $7.6\times 10^{-11}$ \\
  $100$ & $10^{-12}$ & 2.87 & 2.00 & 22.81 & 27.68 & 409.96 & 14.81 & $4.0\times 10^{-14}$ \\
  \bottomrule
\end{tabular}}
\end{table}

\cref{fig:rt2d} reports the wall-clock runtime versus $N$ and the
average throughput at the four target precisions, alongside the
periodic FMM \pfmm{}~\cite{pei_fast_2023} on the same
systems.  The runtime is linear in~$N$ at every precision, and PDMK delivers
roughly $10$- to $26$-fold speedups over \pfmm{} at every
$(N, \varepsilon)$.

\cref{tab:2d} reports the per-step runtime of PDMK and the
head-to-head comparison against \pfmm{} at $N = 10^7$ for the three
aspect ratios $A \in \{1, 10, 100\}$ and the four target
precisions $\varepsilon \in \{10^{-3}, 10^{-6}, 10^{-9}, 10^{-12}\}$.
Three features stand out from the per-step decomposition.  First,
the near-field DMK ($t_{\mathrm{eval}}$) dominates the
per-evaluation cost at every $(A, \varepsilon)$, contributing
$59$--$83\%$ of the total and rising monotonically with
$\varepsilon$ as the PSWF bandwidth grows.  Second, the NUFFT step
is a relatively small and stable contributor at $6$--$12\%$ of the
total across all precisions, reflecting the fact that the 2D
Fourier grid has only modest mode counts at the precisions tested.
Third, $t_{\mathrm{build}}$ (per-evaluation tree rebuild with
image-source placement) is the second-largest contributor at low
precision---up to $34\%$ at $A=1, \varepsilon=10^{-3}$ where the
total budget is small---and falls to $10$--$16\%$ at high precision
where the evaluation cost grows much faster than the build.  The achieved
relative $L^2$ error meets the requested tolerance at every
configuration with $1$--$3$ digits of comfort margin.  The
PDMK-vs-\pfmm{} margin grows with the aspect ratio, peaking near
$A=10$--$100$ at $\varepsilon = 10^{-3}$ ($\approx 25\times$),
and shrinks slightly as $\varepsilon$ tightens because PDMK's
near-field cost grows faster with precision than \pfmm{}'s, which
over-converges by $1$--$2$ digits at high precision.

\subsection{Three dimensions}
\label{sec:results3d}

\begin{figure}[!th]
  \centering
  \includegraphics[width=0.9\linewidth]{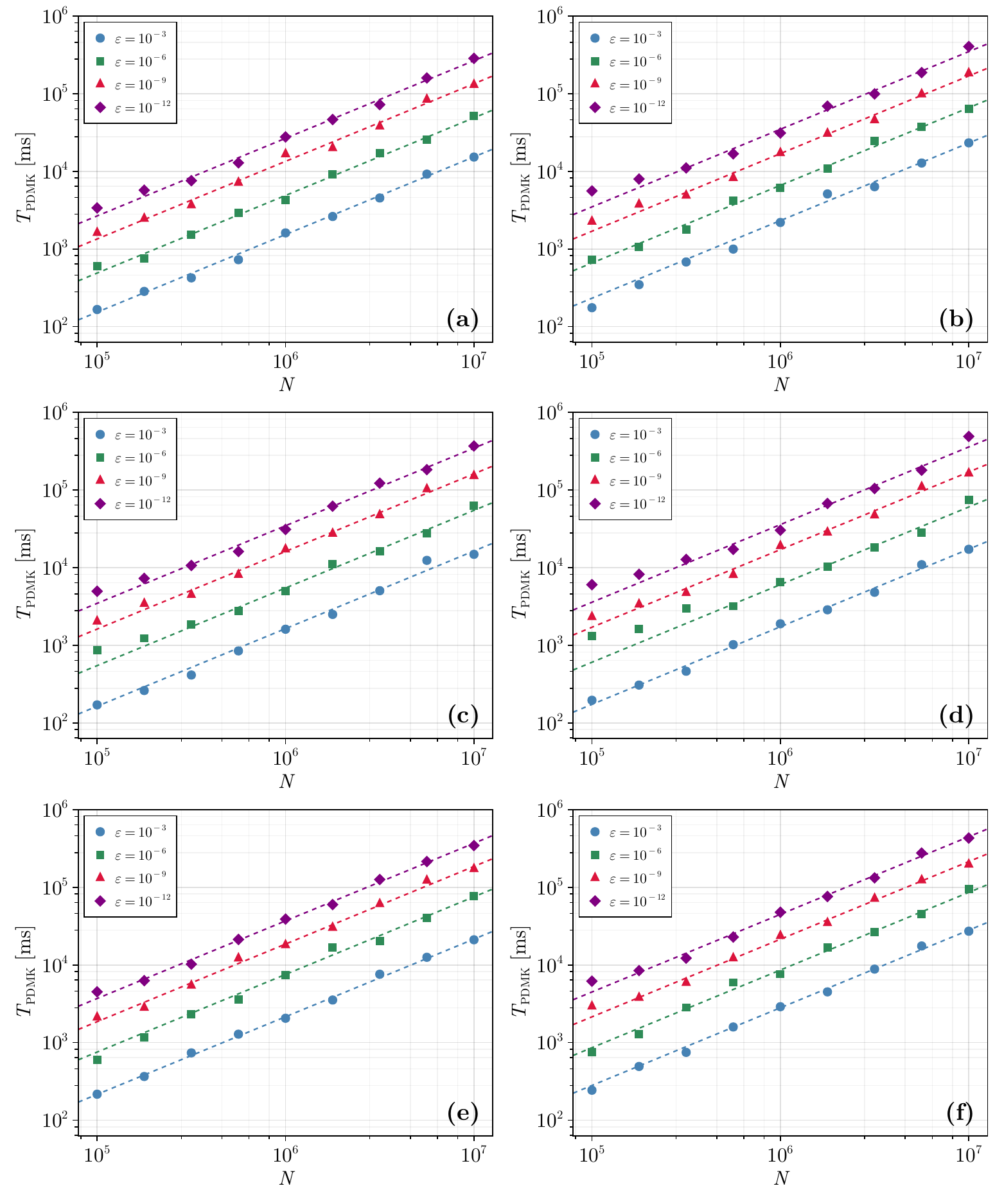}
  \caption{3D runtime versus $N$ in the three scaling regimes; in
    each row the left panel is the rectangular base cell and the
    right panel its triclinic counterpart.  All three rows use the
    same base axes $(1,\,\sqrt{2},\,\sqrt{3})$; rectangular cells
    have $\pi/2$ angles, triclinic $\pi/3$.
    \emph{Top:} isotropic regime~(a); base axes unscaled, so the
    edge ratio stays near $\sqrt{3}$.
    \emph{Middle:} rod regime~(b); the third edge stretched
    $10\times$, giving a long axis roughly $17\times$ the shortest.
    \emph{Bottom:} slab regime~(c); the second and third edges
    each stretched $10\times$, giving two long edges of $14\times$
    and $17\times$ the short edge.  The cell shape within each
    regime is fixed and does not scale with $N$; only $N$ varies
    along each curve.}
  \label{fig:rt3d}
\end{figure}

\begin{table}[!th]
  \centering
  \caption{Per-step PDMK runtime and head-to-head comparison with
    the free-space DMK on the \emph{rectangular} base cell at
    $N = 10^7$, across the three scaling modes and four target
    precisions.  Columns $\varepsilon$, $t_{\mathrm{build}}$,
    $t_{\mathrm{nufft}}$, $t_{\mathrm{eval}}$, $t_{\mathrm{PDMK}}$,
    and $E_{\text{rel}}$ are as defined in \cref{tab:2d}.
    Regime: scaling mode (a), (b), or (c) of \cref{sec:numerics}.
    $t_{\mathrm{DMK}}^0$: total free-space DMK time on the same
    sources.
    ratio: PDMK throughput divided by free-space DMK throughput,
    which (since both run on the same $N$) equals
    $t_{\mathrm{DMK}}^0/t_{\mathrm{PDMK}}$; values above $1$ mean
    PDMK has higher throughput than the free-space baseline.
    Times in seconds; single thread, double precision.}
  \label{tab:rect3d}
  \begin{tabular}{llrrrrrrr}
  \toprule
  Regime & $\varepsilon$
    & $t_{\mathrm{build}}$ & $t_{\mathrm{nufft}}$ & $t_{\mathrm{eval}}$
    & $t_{\mathrm{PDMK}}$ & $t_{\mathrm{DMK}}^0$ & ratio
    & $E_{\text{rel}}$ \\
  \midrule
  (a) iso & $10^{-3}$ & 2.76 & 1.69 & 10.82 & 15.27 & 23.80 & 1.56 & $2.9\times 10^{-3}$ \\
  (a) iso & $10^{-6}$ & 3.83 & 3.12 & 44.82 & 51.77 & 76.27 & 1.47 & $5.7\times 10^{-7}$ \\
  (a) iso & $10^{-9}$ & 4.64 & 7.70 & 118.03 & 130.37 & 185.45 & 1.42 & $7.6\times 10^{-9}$ \\
  (a) iso & $10^{-12}$ & 8.95 & 14.87 & 261.31 & 285.13 & 390.98 & 1.37 & $4.5\times 10^{-13}$ \\
  \midrule
  (b) rod & $10^{-3}$ & 3.02 & 1.71 & 10.12 & 14.84 & 22.93 & 1.54 & $2.4\times 10^{-4}$ \\
  (b) rod & $10^{-6}$ & 4.54 & 5.64 & 52.18 & 62.36 & 75.10 & 1.20 & $6.3\times 10^{-8}$ \\
  (b) rod & $10^{-9}$ & 9.25 & 8.72 & 134.09 & 152.06 & 167.08 & 1.10 & $8.4\times 10^{-10}$ \\
  (b) rod & $10^{-12}$ & 16.50 & 16.80 & 335.59 & 368.90 & 358.95 & 0.97 & $6.8\times 10^{-14}$ \\
  \midrule
  (c) slab & $10^{-3}$ & 5.64 & 1.79 & 13.69 & 21.11 & 20.24 & 0.96 & $1.0\times 10^{-3}$ \\
  (c) slab & $10^{-6}$ & 9.29 & 4.15 & 63.10 & 76.54 & 65.21 & 0.85 & $2.0\times 10^{-7}$ \\
  (c) slab & $10^{-9}$ & 8.11 & 10.01 & 154.16 & 172.28 & 178.10 & 1.03 & $2.0\times 10^{-9}$ \\
  (c) slab & $10^{-12}$ & 15.21 & 18.17 & 311.77 & 345.15 & 310.20 & 0.90 & $1.7\times 10^{-13}$ \\
  \bottomrule
\end{tabular}
\end{table}

\cref{fig:rt3d} shows the wall-clock runtime versus $N$ in the
three 3D scaling regimes, with solid markers at each of the four
target precisions and dashed linear fits.  In each row the left
panel is the rectangular base cell and the right panel its
triclinic counterpart.  All three regimes display linear scaling
in~$N$ at every precision, with the per-row constant growing only
mildly with the regime label.  The triclinic shear inflates the
total runtime by a factor of roughly $1.1$--$1.5$ relative to the
rectangular cell at the same axes, but does not change the slope.

We next decompose the per-step runtime at fixed $N = 10^7$.
\cref{tab:rect3d} reports the breakdown on the rectangular base
cell, alongside the free-space DMK on the same sources for
reference.

\begin{table}[!th]
  \centering
  \caption{Same as \cref{tab:rect3d} but for the
    \emph{triclinic} base cell at $N = 10^7$.  All three lattice
    vectors are unit length and pairwise at angle $\pi/3$ before
    the mode-dependent anisotropic scaling is applied.}
  \label{tab:tri3d}
  \begin{tabular}{llrrrrrrr}
  \toprule
  Regime & $\varepsilon$
    & $t_{\mathrm{build}}$ & $t_{\mathrm{nufft}}$ & $t_{\mathrm{eval}}$
    & $t_{\mathrm{PDMK}}$ & $t_{\mathrm{DMK}}^0$ & ratio
    & $E_{\text{rel}}$ \\
  \midrule
  (a) iso & $10^{-3}$ & 6.18 & 1.66 & 15.43 & 23.27 & 24.53 & 1.05 & $2.7\times 10^{-3}$ \\
  (a) iso & $10^{-6}$ & 7.56 & 3.05 & 53.16 & 63.77 & 77.48 & 1.22 & $4.8\times 10^{-7}$ \\
  (a) iso & $10^{-9}$ & 10.66 & 7.16 & 165.82 & 183.64 & 195.93 & 1.07 & $4.7\times 10^{-9}$ \\
  (a) iso & $10^{-12}$ & 19.98 & 14.69 & 370.15 & 404.83 & 522.08 & 1.29 & $3.7\times 10^{-13}$ \\
  \midrule
  (b) rod & $10^{-3}$ & 3.78 & 1.71 & 11.75 & 17.25 & 30.97 & 1.80 & $3.6\times 10^{-4}$ \\
  (b) rod & $10^{-6}$ & 5.27 & 4.85 & 64.11 & 74.22 & 71.36 & 0.96 & $9.9\times 10^{-8}$ \\
  (b) rod & $10^{-9}$ & 13.72 & 8.60 & 141.64 & 163.96 & 196.23 & 1.20 & $8.2\times 10^{-10}$ \\
  (b) rod & $10^{-12}$ & 14.68 & 21.97 & 452.80 & 489.45 & 310.24 & 0.63 & $1.1\times 10^{-13}$ \\
  \midrule
  (c) slab & $10^{-3}$ & 8.96 & 1.78 & 16.49 & 27.23 & 21.15 & 0.78 & $8.3\times 10^{-4}$ \\
  (c) slab & $10^{-6}$ & 12.82 & 4.27 & 78.59 & 95.68 & 72.98 & 0.76 & $1.5\times 10^{-7}$ \\
  (c) slab & $10^{-9}$ & 13.31 & 9.61 & 174.35 & 197.27 & 223.86 & 1.13 & $1.8\times 10^{-9}$ \\
  (c) slab & $10^{-12}$ & 25.74 & 18.63 & 386.87 & 431.24 & 314.67 & 0.73 & $1.3\times 10^{-13}$ \\
  \bottomrule
\end{tabular}
\end{table}

The near-field DMK time $t_{\mathrm{eval}}$---comprising plane-wave
shifts, upward and downward passes, and leaf-level direct
evaluation of $R_{\mathscr{L}}$---dominates the PDMK runtime at
every configuration, contributing $65$--$92\%$ of
$t_{\mathrm{PDMK}}$ on rectangular cells and a similar share on
triclinic.  The NUFFT step is a uniformly small contributor
at $5$--$12\%$ of $t_{\mathrm{PDMK}}$ across all three regimes and
all four precisions on rect, with the same band on tri; the active
mode count is set by the PSWF bandwidth~$c$ and the cell volume,
which grow only weakly across our regimes.  Tree construction with
image-source placement ($t_{\mathrm{build}}$) is the smallest share
at the precisions of practical interest---$\sim 3$--$8\%$ of
$t_{\mathrm{PDMK}}$ on rectangular cells and $\sim 4$--$13\%$ on
triclinic at $\varepsilon \le 10^{-6}$---but rises to $18$--$33\%$
at $\varepsilon = 10^{-3}$, where the total runtime budget is
small and the image halo and (for triclinic) non-axis-aligned root
grid both contribute.  The relative $L^2$ error meets the
requested tolerance at every configuration with $1$--$3$ digits of
comfort margin.  The corresponding triclinic data are shown
in~\cref{tab:tri3d}.

\begin{figure}[!t]
  \centering
  \includegraphics[width=0.9\linewidth]{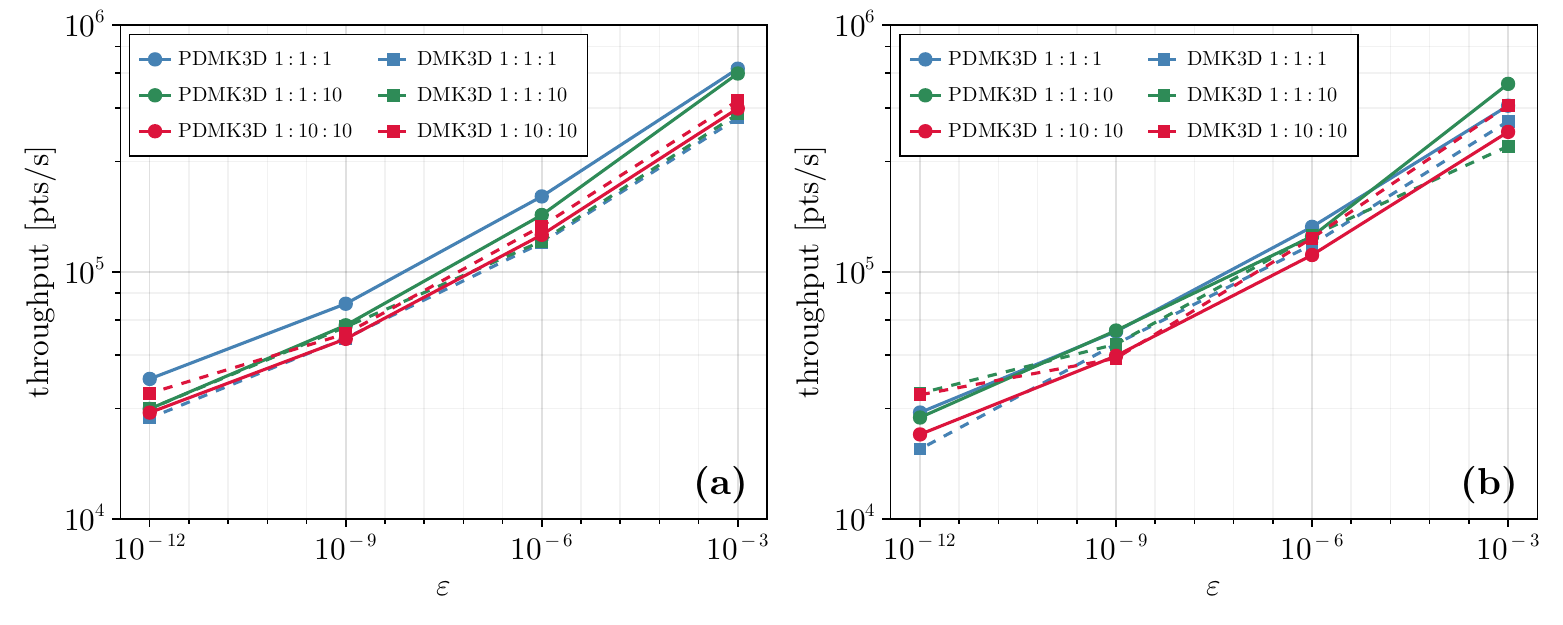}
  \caption{3D average throughput (particles per second)
    for the four target precisions.  (a): rectangular cells (regimes
    (a)--(c)).  (b): triclinic cells with the same axes.}
  \label{fig:tp3d}
\end{figure}

Comparing the PDMK throughput with that of the free-space DMK on
the same sources, PDMK is faster on $9$ of the $12$ rectangular
configurations (throughput ratio $1.03$--$1.56$); slower performance
occurs only in the (c)~slab regime, where the ratio in the interval
$0.85$--$0.96$.  On triclinic cells,the ratios spread more widely,
$0.63$--$1.80$, with the (b)~rod regime favoring PDMK at low
precision but the (c)~slab regime relative performance
falling to $0.7$--$0.8$ as the
non-axis-aligned root grid inflates the image halo.  \cref{fig:tp3d}
summarizes the average throughput across the three regimes and
four precisions.

% ---------------------------------------------------------
\section{Conclusions}
\label{sec:conclusions}
% ---------------------------------------------------------

We have extended the DMK
framework to periodic boundary conditions in arbitrary unit cells
(oblique in 2D, triclinic in 3D).  The single-cube root of
free-space DMK is replaced by a rectangular grid covering the unit
cell plus a one-cube image halo, on which the compactly supported
components of the kernel splitting are evaluated by a free-space
DMK; the long-range periodic kernel $W_0^{(p)}$ is represented by a single
Fourier series on the reciprocal lattice and evaluated by the
five-step procedure of fast Ewald summation.  Quasi-periodic
(Bloch) boundary conditions are supported through three localized
modifications.

Compared with the periodic FMM of~\cite{pei_fast_2023}, our scheme
uses a single Fourier series rather than $2d$ directional plane-wave
expansions, and the root-cube size $\rc$ is a free parameter that
adapts to source nonuniformity and unifies the FMM and fast Ewald
summation at the level of implementation.  At fixed cell shape, the
algorithm has $O(N)$ complexity.

Numerical experiments on highly nonuniform 2D and 3D source
distributions confirm the analysis.  In two dimensions on oblique
cells with aspect ratios up to $100$, our algorithm delivers
$10$- to $26$-fold speedups over \pfmm{}~\cite{pei_fast_2023}
across all particle counts and target precisions tested.  In three
dimensions, on both rectangular and triclinic cells with edge-length
ratios up to roughly $17$, PDMK throughput stays within
$0.6$--$1.8\times$ of the underlying free-space DMK on the same
sources.

Several extensions are natural.  The DMK framework is largely
kernel-agnostic, so the construction extends to other Green's
functions with only modest changes to the top-level periodic
kernel.  Mixed boundary conditions---in particular slab
geometries (periodic in two directions, free in the third)---are
common in surface and interface problems and can be handled by
adapting the Fourier series to the partially periodic lattice
(see \cite{krantz2026,shamshirgar_fast_2021} and \cite{jiang2025nm}).
Distributed-memory parallelization (MPI) and GPU
acceleration of the tree, plane-wave, and NUFFT passes would
substantially extend the range of accessible system sizes.

\bibliography{pbcdmk}

\end{document}